\documentclass[12pt,reqno]{amsart}
\usepackage[letterpaper,margin=1.15in]{geometry}
\usepackage{amssymb,mathtools}
\usepackage{microtype}
\usepackage{xcolor}
\usepackage{graphicx,accsupp}
\usepackage[section]{placeins}
\usepackage{hyperref}

\definecolor{linkblue}{RGB}{0,85,170}
\hypersetup{
  colorlinks=true,
  linkcolor=linkblue,
  citecolor=linkblue,
  urlcolor=linkblue,
  pdftitle={Brunn-Minkowski inequalities for Riesz capacity in all dimensions},
  pdfauthor={Carrie Clark and Richard S. Laugesen},
  pdfsubject={Riesz capacity in higher dimensions, weighted extension and equality},
  pdfkeywords={Riesz capacity, logarithmic capacity, Brunn-Minkowski inequality,
  fractional Laplacian, optimal transport}
}

\newtheorem{theorem}{Theorem}
\newtheorem{lemma}[theorem]{Lemma}
\newtheorem{corollary}[theorem]{Corollary}
\newtheorem{conjecture}[theorem]{Conjecture}
\theoremstyle{definition}
\newtheorem{example}[theorem]{Example}
\newtheorem{problem}[theorem]{Open Problem}
\theoremstyle{plain}
\newcommand{\R}{\mathbb R}
\newcommand{\PP}{\mathcal P}
\DeclareMathOperator{\RCap}{Cap}
\DeclareMathOperator{\diam}{diam}
\DeclareMathOperator{\supp}{supp}
\newcommand{\doi}[1]{\href{https://doi.org/#1}{doi:\nolinkurl{#1}}}

\title[Brunn--Minkowski for Riesz capacity]{Brunn--Minkowski inequalities for Riesz
capacity in all dimensions}
\author{Carrie Clark}
\address{Kitchener, Ontario, Canada}
\email{carrieclark435@gmail.com}
\author{Richard S. Laugesen}
\address{University of Illinois, Urbana IL 61801, USA}
\email{Laugesen@illinois.edu}
\date{\today}

\newcommand{\La}{\mathcal L}
\newcommand{\grad}{\nabla}
\newcommand{\flatset}[1]{\widehat{#1}}
\newcommand{\dd}{\,d}

\begin{document}

\begin{abstract}
We prove the Brunn--Minkowski inequality for Riesz $p$-capacity on compact
subsets of $\R^n$ when $n-1<p<n$, extending the known endpoint $p=n-1$.
For sets of positive capacity, equality holds if and only if they are
convex and positively homothetic. The inequality remains conjectural for
$n-2<p<n-1$, when $n\ge2$. Known isodiametric counterexamples rule out a
general extension to negative $p$, even for convex planar sets.

\end{abstract}

\maketitle

\section{\bf Introduction and main results}\label{sec:introduction}

Is Riesz capacity superlinear on compact sets? Such a Brunn--Minkowski inequality indeed
holds in one dimension for Riesz exponents less than $1$, as we recently proved
\cite{CLline}. On $n$-dimensional sets, the current paper establishes capacitary
Brunn--Minkowski for exponents between $n-1$ and $n$. Equality cases are characterized.
Known isodiametric counterexamples constrain the possible extension to negative
exponents.

\subsection*{Riesz capacity}

For a nonempty compact set $K\subset\R^n$, write $\PP(K)$ for its
Borel probability measures and put
\[
  \RCap_p(K)=V_p(K)^{-1/p} , \qquad V_p(K)=\min_{\mu\in\PP(K)}
  \iint_{K\times K}|x-y|^{-p}\dd\mu(x)\dd\mu(y) ,
\]
when $p>0$, with the convention that $\infty^{-1/p}=0$. For $p<0$ the minimum is
replaced by a maximum. The double integral is called the
$p$-energy of $\mu$, while $V_p(K)$ is the $p$-energy of $K$ and $\RCap_p(K)$ is its
Riesz $p$-capacity. At $p=0$
we use logarithmic capacity,
\[
  \RCap_0(K)=e^{-V_0(K)} , \qquad V_0(K)=\min_{\mu\in\PP(K)}
  \iint_{K\times K}\log\frac1{|x-y|}\dd\mu(x)\dd\mu(y) .
\]
We restrict to $p<n$, since $\RCap_p(K)=0$ for every
compact $K\subset\R^n$ when $p\ge n$
\cite[Theorem 4.3.1]{BHS19}; see also
\cite[Theorem 1.1(d)]{CLproperties}.

Measures attaining these extrema are called equilibrium measures.
For Riesz potential theory, see Landkof \cite{L72} and
Borodachov--Hardin--Saff \cite[Chapter 4]{BHS19}.
Our capacity normalization differs from theirs: the exponent $-1/p$ in our definition
makes capacity homogeneous of degree one. Capacity is invariant
under translations and orthogonal maps.

The Minkowski sum of nonempty sets $K,L$ is
\[
  K+L=\{x+y:x\in K,\ y\in L\} .
\]

Brunn--Minkowski is known in one dimension:

\begin{theorem}[Capacitary Brunn--Minkowski on the line; Clark--Laugesen
\protect{\cite[Thm.\ 1]{CLline}}]\label{thm:line}
If $p<1$ and $K,L\subset\R$ are nonempty compact sets, then
\[
  \RCap_p(K+L)\ge\RCap_p(K)+\RCap_p(L).
\]
When $p\le-1$, equality holds for all such $K,L$.
\end{theorem}

The one dimensional quantile argument in \cite{CLline} relies on order
preservation. Section~\ref{sec:transport} below extends that trial measure
construction to higher dimensions through transport interpolation, obtaining a weaker
inequality.

For general compact sets, a different method yields Brunn--Minkowski in all dimensions
provided the Riesz exponent lies within $1$ of the dimension $n$.

\begin{theorem}[Capacitary Brunn--Minkowski above $n-1$]\label{thm:upper-band}
Let $n\ge1$. If $n-1\le p<n$ and $K,L\subset\R^n$ are nonempty compact
sets, then
\begin{equation}\label{eq:sumBM}
  \RCap_p(K+L)\ge\RCap_p(K)+\RCap_p(L).
\end{equation}
\end{theorem}

By homogeneity, inequality \eqref{eq:sumBM} for all $K,L$ is equivalent to the concavity
inequality
\begin{equation}\label{eq:higher-dimensional-BM}
  \RCap_p\bigl((1-\lambda)K+\lambda L\bigr)
  \ge(1-\lambda)\RCap_p(K)+\lambda\RCap_p(L),\qquad 0\le\lambda\le1.
\end{equation}

The proof of the theorem in Sections~\ref{sec:extension} and \ref{sec:upper-proof}
applies the Caffarelli--Silvestre extension \cite{CS07} to equilibrium
potentials for $(-\Delta)^{(n-p)/2}$, yielding a local weighted equation
in one additional variable.

The new range covered by the theorem when $n\ge2$ is $n-1<p<n$.
Akman \cite[Theorems A(a) and B(e)]{A26} proved the Newtonian case $p=n-2$
for arbitrary compact sets, including logarithmic capacity when $n=2$.
Embedding $\R^n$ into $\R^{n+1}$ then handles the endpoint case $p=n-1$.
When $n=1$, Theorem~\ref{thm:line} already covers all $p<1$, and $p=0$ follows also from
Akman's logarithmic theorem.
Section~\ref{sec:literature} describes these and other predecessors.

\subsection*{The volume limiting case}
For each compact $M\subset\R^n$, the capacity-to-volume limit is
\[
  \lim_{p\uparrow n}
  \left(\frac{|\mathbb S^{n-1}|}{n-p}\right)^{\!\!1/p}\RCap_p(M)
  =m_n(M)^{1/n}
\]
by \cite[Corollary 1.2]{CLproperties}, where $m_n$ is Lebesgue measure.
Multiplying \eqref{eq:sumBM} by the quantity
$(|\mathbb S^{n-1}|/(n-p))^{1/p}$ and letting $p\uparrow n$ recovers
Brunn--Minkowski for volume:
\[
  m_n(K+L)^{1/n}\ge m_n(K)^{1/n}+m_n(L)^{1/n}.
\]

\subsection*{Remaining conjectural range} To fill in the gap between Akman's result at
$p=n-2$ and \autoref{thm:upper-band} for $n-1 \leq p < n$, one raises:

\begin{conjecture}[Capacitary Brunn--Minkowski below $n-1$]\label{conj:higher-dimensional}
Let $n\ge2$. If $n-2<p<n-1$ and $K,L\subset\R^n$ are nonempty compact
sets, then \eqref{eq:sumBM} holds.
\end{conjecture}

For convex sets, this range was suggested in
\cite[Section 10, Remark 1]{CLratios}. Novaga--Ruffini
\cite[Conjecture 3.3]{NR15} had earlier proposed $0<p<n$.

\subsection*{Equality}

The following theorem characterizes equality for summands of positive
capacity. Here positive homothety means $L=cK+b$ for some $c>0$ and
$b\in\R^n$.

\begin{theorem}[Equality iff convexity and positive homothety]\label{thm:upper-equality}
Let $n\ge1$ and $n-1\le p<n$. If $K,L\subset\R^n$ are compact sets
with positive $p$-capacities, then equality holds in \eqref{eq:sumBM}
if and only if $K,L$ are convex and positively homothetic.
\end{theorem}

This equality characterization also applies to
\eqref{eq:higher-dimensional-BM} for $0<\lambda<1$.

The proof in Section~\ref{sec:upper-equality} adapts methods of
Caffarelli--Jerison--Lieb \cite{CJL96}, Colesanti--Salani \cite{CS03}
and Akman \cite{A26} to the weighted extension equation.
The equality characterization at $p=n-1$ follows from Akman
\cite[Theorems A(d) and B(e)]{A26} by embedding in $\R^{n+1}$,
as explained in Section~\ref{sec:literature}. For $n=1$, the equality
characterization throughout the range $0\le p<1$ is also proved in
\cite[Theorem 2]{CLline}.

\subsection*{Mean width and diameter}

Write $B_n$ for the closed unit ball in $\R^n$ and $w(K)$ for the mean directional
width of $K$: the spherical average of the distance between parallel
supporting hyperplanes of its convex hull. Note $w(B_n)=2$.

\begin{corollary}[Capacity, mean width and diameter]\label{cor:width-diameter}
If $n\ge1$, $n-1\le p<n$ and $K\subset\R^n$ is compact with positive
diameter, then
\begin{equation}\label{eq:width-diameter}
  \frac{\RCap_p(K)}{\RCap_p(B_n)}
  \le\frac{w(K)}2\le\frac{\diam K}2.
\end{equation}
Equality in the first inequality holds if and only if $K$ is a closed
ball. The same condition characterizes equality between the first and
last quantities. Equality in the second inequality alone means that
the convex hull of $K$ has constant width.
\end{corollary}

Section~\ref{sec:upper-equality} provides a direct
proof based on capacitary Brunn--Minkowski.

For $n=1$, the corollary is \cite[Proposition 3]{CLratios}:
$w(K)=\diam K$, and the capacity bound follows for all $p<1$ by containing $K$
in its convex hull, which is an interval.

For $n\ge2$, on the wider range $n-2\le p<n$ the inequalities in the corollary can be
obtained from results in the literature that do not use Brunn--Minkowski. At the lower
endpoint $p=n-2$ with $n \geq 3$,
the Newtonian bounds go back to Szeg\H{o}
\cite[S\"atze II und III]{S31}, whose argument extends to all $n\ge3$;
see also Borell \cite[Example 7.5]{B84}. For the planar logarithmic
bounds with $p=0, n=2$, see P\'olya--Szeg\H{o} \cite[p.~17, entry 12 and p.~19, entry
31]{PS51}
and Borell \cite[Example 7.4]{B84}. (P\'{o}lya and Szeg\H{o}'s outer conformal radius is our logarithmic capacity, and their perimeter $L$ equals $\pi w$ for planar convex bodies.) On the full range $n-2<p<n$,
Su--Sun--Xiong \cite[Theorem 1.1 and (2.1)]{SSX26} obtain the first inequality
(``extremal ratio of capacities'') in
\[
  \frac{\RCap_p(K)}{\RCap_p(B_n)}
  \le\frac{\RCap_{n-2}(K)}{\RCap_{n-2}(B_n)}
  \le\frac{w(K)}2
\]
while the second is the Newtonian or logarithmic mean width bound just
cited. Hence both bounds hold throughout the larger range.
The equality characterization fails
at $p=n-2$, however, because a sphere has the same width, diameter and Newtonian
capacity as a solid ball.
Incidentally, the Brunn--Minkowski deduction of the isodiametric inequality when
$p=n-2,n-1$ was noted in \cite[Theorem 13 and Section 10]{CLratios}. Novaga--Ruffini
\cite[Corollary 3.2]{NR15} previously gave the mean width bound at
$p=n-1$.

\subsection*{Failure for negative \texorpdfstring{$p$}{p}}

The square in Example~\ref{ex:square} gives a simple planar counterexample
for $p\le-1$. More generally, the isodiametric counterexamples of
Clark--Laugesen \cite[Theorem 7]{CLratios} imply failure even for convex
planar summands when $p<q_*=-0.856\ldots$;
see Corollary~\ref{cor:isodiametric-failure}.
Open Problem~\ref{prob:planar-negative} asks whether any negative exponent
admits a planar Brunn--Minkowski inequality.

\subsection*{Failure for positive \texorpdfstring{$p$}{p}}
For each fixed $p>0$, Burchard--Choksi--Hess-Childs
\cite[Theorem 4]{BCH20} show that balls fail to maximize $\RCap_p$ at
fixed diameter when $n$ is sufficiently large. Brunn--Minkowski would
imply this isodiametric inequality by taking $L=-K$, and so it too must
fail. These counterexamples lie below the conjectured range $p>n-2$.

\subsection*{Lower bounds from couplings}

Theorem~\ref{thm:transport} uses optimal transport to show that
\[
  \RCap_p(K+L)\ge
  \sqrt{\RCap_p(K)^2+\RCap_p(L)^2} , \qquad p \in [-2,n) ,
\]
for nonempty compact sets $K,L\subset\R^n$. Although weaker than
Brunn--Minkowski, this bound does cover exponents where
Brunn--Minkowski is either open or known to fail.
An independent coupling gives the companion bound
\eqref{eq:product-bound} for $p<-2$.

\subsection*{Structure of the paper}

Section~\ref{sec:literature} reviews prior work.
Section~\ref{sec:extension} develops the weighted extension of equilibrium
potentials and continuity under outer approximation, which
Section~\ref{sec:upper-proof} uses to prove Theorem~\ref{thm:upper-band}
when $n-1\le p<n$. Section~\ref{sec:upper-equality} proves the
convexity and positive homothety characterization in
Theorem~\ref{thm:upper-equality}. It also proves the mean width and isodiametric bounds in
Corollary~\ref{cor:width-diameter}.

Section~\ref{sec:negative-planar-proof} gives the square counterexample
in Example~\ref{ex:square} and explains how known isodiametric
counterexamples imply failure for further negative exponents.
Section~\ref{sec:transport} proves the coupling bounds in
Theorem~\ref{thm:transport}.
Appendix~\ref{app:level-set-foliation} expands the supporting point and
foliation arguments from Section~\ref{sec:upper-equality}.

\section{\bf Earlier capacitary inequalities}\label{sec:literature}

We review known inequalities and equality cases, distinguishing
Riesz capacities from variational capacities.

\subsection*{Newtonian and logarithmic capacities}

Pommerenke \cite[Satz 2]{P59} proved the Brunn--Minkowski inequality
for logarithmic capacity of connected planar compact sets. His book
\cite[Corollary 11.1 and Theorem 11.2(b)]{P75} can be used to provide an alternative proof.
Borell proved the Brunn--Minkowski inequality for Newtonian capacity
($p=n-2$) on convex bodies in $\R^n$, $n\ge3$ \cite{B83},
and for logarithmic capacity ($p=0$) on convex planar sets
\cite[Example 7.3]{B84}. For bodies with
nonempty interiors and $0<\lambda<1$, equality means positive homothety,
by Caffarelli--Jerison--Lieb \cite[Theorem B]{CJL96} in the Newtonian
case and Colesanti--Cuoghi \cite[Theorem 3.1]{CC05} in the planar logarithmic case.

For planar logarithmic capacity, Borell
\cite[Theorem 7.2 and Example 7.3]{B84} uses the same maximum principle
method as in his Newtonian proof. He first compares equilibrium potentials
in expanding disks and then passes to exterior Green functions with pole
at infinity. Comparing the constant terms in their logarithmic expansions
at infinity gives the bound by the minimum of the summand capacities;
homogeneity then gives Brunn--Minkowski. 

Novaga--Ruffini \cite[Theorem 1.1(ii)]{NR15} treated $p=n-1$ for
convex bodies through harmonic extension. Alternatively, embedding into $\R^{n+1}$
reduces this case to Borell's theorem \cite[p.~1064]{MNR18}, as remarked around
\cite[Theorem 20]{CLratios}. Muratov--Novaga--Ruffini
\cite[Lemma 3.7]{MNR18} characterized equality by positive homothety
in the planar case with $p=1$.

Akman \cite[Theorems A(a), A(d) and B(e)]{A26} was the first to prove the
Brunn--Minkowski inequality on arbitrary compact
sets, for Newtonian and planar logarithmic capacity, that is, $p=n-2$ for $n \geq 2$.
For positive
capacities and $0<\lambda<1$, he showed equality exactly when the filled sets
are convex and positively homothetic. Again, these results extend to $p=n-1$ by
embedding $K,L\subset\R^n$ into $\R^{n+1}$. The embedded sets have no bounded
complementary components, and so equality there means that $K,L$
themselves are convex and positively homothetic.

\subsection*{Variational capacities}

Colesanti--Salani \cite[Theorem 1]{CS03} proved Brunn--Minkowski for
variational $r$-capacity when $1<r<n$, with positive homothety
characterizing equality for convex bodies with nonempty interiors.
Colesanti--Cuoghi \cite{CC05} treated the analogous $n$-Laplace
logarithmic capacity and its equality cases. Their capacity agrees with
$\RCap_0$ when $n=2$, but differs from our logarithmic capacity when $n>2$.
For homogeneous nonlinear
elliptic operators, Akman--Gong--Hineman--Lewis--Vogel
\cite[Theorem A]{AGHLV22} proved the inequality for convex compact
sets of positive capacity, without requiring nonempty interior. Akman--Lewis--Saari--Vogel
\cite[Theorem A]{ALSV21} developed Green function capacities for
$r\ge n$. Both works characterize equality by positive homothety under
operator hypotheses satisfied by the $r$-Laplacian, for summands of
positive capacity. Akman's compact-set theorem extends the
$r$-Laplace results beyond convex sets.

The variational exponent $r$ differs from the Riesz kernel exponent $p$.
The families coincide in the Dirichlet integral case $r=2$, corresponding
to Newtonian capacity when $n\ge3$ and logarithmic capacity when $n=2$.

\subsection*{Capacity ratios}

Clark--Laugesen \cite{CLratios} studied maximization of
$\RCap_q(K)/\RCap_p(K)$ among compact subsets of $\R^n$,
with balls, regular simplex vertex sets and two-point sets as
proved or conjectured maximizers in different parameter regions.
A technical report of Liu \cite[Corollary 1.3]{L26} proves
two-point maximality for $q<p<0$ and regular simplex maximality
for $p<q\le-2$, in every dimension.
Su--Sun--Xiong \cite[Theorems 1.1 and 1.2]{SSX26} prove ball
maximality for $\RCap_q(K)/\RCap_{n-2}(K)$ when $n\ge2$ and
$n-2<q<n$, and for $\RCap_{n-2}(K)/\RCap_0(K)$ when $n\ge4$.
The former includes the planar P\'olya--Szeg\H{o} capacity ratio
conjecture when $(n,q)=(2,1)$.

\section{\bf Extended equilibrium potentials: the Caffarelli--Silvestre
extension}\label{sec:extension}

We apply the Caffarelli--Silvestre extension \cite{CS07} to equilibrium
potentials and use classical Riesz potential theory to control their
boundary values. Parallel neighborhoods provide the regularity needed
to approximate arbitrary compact sets.

For the extension argument in this section, consider
\[
  n\ge1,\qquad p>0,\qquad n-2<p<n .
\]
Write $X=(x,t)\in\R^{n+1}$ and
\[
  \flatset K=K\times\{0\}
\]
for the flat copy of a compact set $K\subset\R^n$. Suppose $K$ has positive capacity,
that is, finite energy $0<V_p(K)<\infty$. Take an equilibrium measure $\mu_K$ and define
\begin{equation}\label{eq:extended-potential}
  U_K(x,t)=\frac1{V_p(K)}
  \int_K (|x-y|^2+t^2)^{-p/2}\dd\mu_K(y).
\end{equation}
This is the normalized $p$-equilibrium potential of the flat set
$\flatset K$, equal to $1$ there except perhaps on a set of zero $p$-capacity.
The kernel here is a constant multiple of the weighted fundamental solution
in \cite[Section 2.2]{CS07}.

Direct differentiation gives the equation
\begin{equation}\label{eq:weighted-equation}
  \La U_K:=\Delta_{x,t}U_K+\frac a t(U_K)_t=0,\qquad t \neq 0 ,
\end{equation}
where
\[
  a=p-n+1 \in(-1,1) .
\]
For $t\ne0$, the equivalent divergence form of the equation is
$\operatorname{div}(|t|^a\grad U_K)=0$.

\subsection*{Horizontal regularization of the set}

Parallel neighborhoods allow us to approximate arbitrary compact sets,
following a strategy used by Akman
\cite[Lemmas 4.9 and 4.10]{A26} for variational capacity.
For nonempty compact $E$ and $r>0$, the parallel neighborhood
$K=E+rB_n$ contains a ball and hence has positive $p$-capacity.

The next lemma gives continuity up to the regularized flat set and identifies the
capacity from the potential's decay at infinity.

\begin{lemma}[Potentials of parallel neighborhoods]\label{lem:extended-facts}
Suppose $n \geq 1$, $p>0$ and $n-2<p<n$. Let $E\subset\R^n$ be nonempty and compact and
let $r>0$. If $K=E+rB_n$
then $U_K$ is continuous on $\R^{n+1}$, smooth off $\flatset K$ and even
in $t$. Moreover,
\[
  0<U_K\le1,\qquad \{U_K=1\}=\flatset K,
\]
and
\begin{equation}\label{eq:potential-infinity}
  \lim_{|X|\to\infty}|X|^p \, U_K(X)
  = \frac{1}{V_p(K)}=\RCap_p(K)^p .
\end{equation}
Further, $U_K$ has no local maximum outside $\flatset K$. Equation
\eqref{eq:weighted-equation} holds classically when $t\ne0$, and the divergence-form
equation $\operatorname{div}(|t|^a\grad U_K)=0$ holds weakly
in $\R^{n+1}\setminus\flatset K$.
\end{lemma}

\begin{proof}
Regard $\mu_K$ as a measure in $\R^{n+1}$. Compactness of $K$
gives the limit \eqref{eq:potential-infinity} directly from
\eqref{eq:extended-potential}. The Frostman inequalities and maximum
principle \cite[Theorems 4.2.14, 4.2.15 and 4.5.11]{BHS19}, applied
in $\R^n$, give $0<U_K(x,0)\le1$ everywhere and $U_K(x,0)=1$ on
$K$ outside an exceptional set of zero $p$-capacity. The integral definition
\eqref{eq:extended-potential} gives
$0<U_K(x,t)\le U_K(x,0)\le1$ for all $(x,t)$.

To remove the exceptional set, fix $x\in K$ and choose a closed ball
$D\subset K$ of radius $r$ containing $x$. Since $p<n$, Lebesgue
measure on $D$ has finite $p$-energy and therefore does not charge the
exceptional set, by \cite[Lemma 4.2.13]{BHS19}. Hence
$U_K(\cdot,0)=1$ Lebesgue-almost everywhere on $D$.
For $0<\varepsilon<r$, put $D_\varepsilon=D\cap B(x,\varepsilon)$.
The set $D_\varepsilon$ contains a ball of radius $\varepsilon/2$,
and so $|D_\varepsilon|\ge c\varepsilon^n$, even if $x\in\partial D$.

We claim that
\begin{equation}\label{eq:average}
  \frac1{|D_\varepsilon|}\int_{D_\varepsilon}|z-y|^{-p}\dd z
  \le C|x-y|^{-p},\qquad
  y\in\mathbb R^n\setminus\{x\},
\end{equation}
with $C$ independent of $\varepsilon$ and $y$, as we show in two cases.
If $|x-y|\ge2\varepsilon$, then $|z-y|\ge|x-y|/2$ throughout
$D_\varepsilon$. If $|x-y|<2\varepsilon$, then
$D_\varepsilon\subset B(y,3\varepsilon)$, and integration over that
ball bounds the average by $C\varepsilon^{-p}\le C'|x-y|^{-p}$.

For each $y\ne x$, continuity of the kernel near $x$ makes this average
converge to $|x-y|^{-p}$ as $\varepsilon\downarrow0$.
The bound $U_K(x,0)\le1$ implies that $\mu_K$ has no atom at $x$ and
that the dominating function in \eqref{eq:average} is $\mu_K$-integrable.
Fubini's theorem and dominated convergence therefore give
\[
  1=\frac1{|D_\varepsilon|}\int_{D_\varepsilon}U_K(z,0)\dd z
  \longrightarrow
  \frac1{V_p(K)}\int_K|x-y|^{-p}\dd\mu_K(y)=U_K(x,0).
\]
Consequently $U_K=1$ on $\flatset K$.

Lower semicontinuity in $\R^{n+1}$ and the bound $U_K\le1$ now give continuity at
$\flatset K$. Smoothness of $U_K$ is immediate from the definition, away from $\flatset
K$, as is
evenness in $t$. The $t$-derivative is
\begin{equation}\label{eq:vertical-derivative}
  (U_K)_t(x,t)=-\frac{pt}{V_p(K)}
  \int_K (|x-y|^2+t^2)^{-(p+2)/2}\dd\mu_K(y)
\end{equation}
for $t\ne0$. At each point $(x,0)\notin\flatset K$,
\begin{equation}\label{eq:horizontal-laplacian}
  \Delta_xU_K(x,0)=\frac{p(p-n+2)}{V_p(K)}
  \int_K |x-y|^{-p-2}\dd\mu_K(y)>0.
\end{equation}
Here the strict inequality $p>n-2$ is used.
The $t$-derivative in \eqref{eq:vertical-derivative} excludes local maxima when $t\ne0$,
while
the horizontal Laplacian in \eqref{eq:horizontal-laplacian} excludes them when $t=0$,
outside
$\flatset K$. Hence $\{U_K=1\}=\flatset K$.

For $\psi\in C_c^\infty(\R^{n+1}\setminus\flatset K)$,
integration by parts over $|t|>\varepsilon$, using
$\operatorname{div}(|t|^a\grad U_K)=0$ there, gives
\[
  \begin{split}
    \left|\int_{|t|>\varepsilon}|t|^a\grad U_K\cdot\grad\psi\dd X\right|
    &\le\varepsilon^a\sum_{\pm}\int_{\R^n}
    \left|\psi(x,\pm\varepsilon)(U_K)_t(x,\pm\varepsilon)\right|\dd x\\
    &\le C_\psi\varepsilon^{a+1}\longrightarrow0
  \end{split}
\]
as $\varepsilon \to 0$, since $a+1=p-n+2>0$ and \eqref{eq:vertical-derivative} gives
$(U_K)_t=O(t)$ uniformly on the
support of $\psi$ near $t=0$. Since $a>-1$ and $\grad U_K$ is bounded on this support,
dominated
convergence gives the weak equation $\int_{\R^{n+1}}|t|^a\grad U_K\cdot\grad\psi\dd X=0$.
\end{proof}

\subsection*{Continuity under outer approximation of the set}

The following continuity result enables us to pass back to the original compact
set. The logarithmic case $p=0$ is allowed here too.

\begin{lemma}[Decreasing compact sets]\label{lem:outer-continuity}
Let $n \geq 1$ and $p \geq 0$. If the $K_j$ are nonempty compact sets in $\R^n$ with
$K_j \supset K_{j+1}$ and $K=\bigcap_jK_j$, then
$\RCap_p(K_j)\downarrow\RCap_p(K)$ as $j \to \infty$.
If $0<p<n$ and $\RCap_p(K)>0$, then the equilibrium measures
$\mu_j$ of $K_j$ converge weak-$*$ to $\mu_K$ and their extended
potentials $U_j$ converge locally uniformly to $U_K$ away from $\flatset K$.
\end{lemma}

\begin{proof}
By inclusion monotonicity, the energies $V_p(K_j)$ increase to a limit $\ell\le V_p(K)$.
If $\ell=\infty$, the conclusion is immediate. Suppose $\ell<\infty$ and choose a
weak-$*$ convergent subsequence of equilibrium measures in
$\PP(K_1)$. Its limit is supported on $K$, and the Principle of
Descent \cite[Theorem 4.1.2]{BHS19} (applicable also to the logarithmic kernel) implies
that its energy is at most $\ell$. Hence $V_p(K)\le\ell$, proving $V_p(K_j)\uparrow
V_p(K)$ and hence $\RCap_p(K_j)\downarrow\RCap_p(K)$.

If $0<p<n$ and $\RCap_p(K)>0$, then every weak-$*$ cluster point of the equilibrium
measures $\mu_j$ is an
equilibrium measure on $K$, by the same descent argument.
Uniqueness of the equilibrium measure on $K$ \cite[Theorem 4.4.5]{BHS19} forces convergence of
$\mu_j$ to $\mu_K$.
For every compact set away from $\flatset K$, the sets $\flatset K_j$
are uniformly separated from it for large $j$. Uniform continuity of
the kernel there, together with $V_p(K_j)\to V_p(K)>0$, gives uniform
convergence of $U_j$ to $U_K$.
\end{proof}

\section{\bf Proof of Theorem~\ref{thm:upper-band} --- Brunn--Minkowski for
\texorpdfstring{$n-1\le p<n$}{n-1 <= p < n}}\label{sec:upper-proof}

We prove the theorem for $p>0$ and $n-1\le p<n$. In that range, one has
\[
  0\le a<1 .
\]
The remaining case $(n,p)=(1,0)$ follows from Theorem~\ref{thm:line}.
The endpoint $p=n-1$ is also known from Akman's results by embedding
in $\R^{n+1}$, as explained in Section~\ref{sec:literature}.

Our proof builds on the convolution comparison of Borell
\cite{B83}, Colesanti--Salani \cite{CS03} and Akman
\cite[Section 4]{A26}. Akman's argument takes place in the original ambient space.
Here we adapt the comparison method to the Caffarelli--Silvestre extension
of the Riesz potentials in $\R^{n+1}$, governed by \eqref{eq:weighted-equation}.
Novaga--Ruffini \cite{NR15} combined this comparison with
harmonic extension at $p=n-1$.
The weighted equation requires additional arguments concerning
optimizing gradients, the drift term and passage across $t=0$,
which are provided in Steps~1, 3 and~4 below.

\subsection*{Step 1 --- supremum--minimum convolution and optimizing pairs}

Let $K,L\subset\R^n$ be nonempty compact sets of positive
capacity, fix $0<\lambda<1$ and put
\[
  M=(1-\lambda)K+\lambda L.
\]
Write $U_0=U_K$ and $U_1=U_L$ for their extended potentials, as in
Section~\ref{sec:extension}, and define their supremum--minimum convolution by
\begin{equation}\label{eq:sup-convolution}
  U_\lambda(X)=\sup_{X=(1-\lambda) Y+\lambda Z}
  \min\{U_0(Y),U_1(Z)\}
\end{equation}
where the supremum is taken over all $Y,Z \in \R^{n+1}$ such that $(1-\lambda) Y+\lambda
Z=X$.

First assume that $K=A+rB_n$ and $L=C+rB_n$, where $A,C\subset\R^n$
are nonempty compact sets, $r>0$ and $B_n$ is the closed unit ball in $\R^n$.
Then $M=(1-\lambda)A+\lambda C+rB_n$.

The attainment, continuity and decay conclusions in the following
lemma follow from a standard compactness argument; see
\cite[Lemma 4.1]{A26}. We also establish
local Lipschitz regularity and exclude critical optimizing inputs at
points of nonzero height, as will be needed for the later weighted comparison.

\begin{lemma}[Optimizers and local Lipschitz regularity]\label{lem:optimizers}
The supremum in \eqref{eq:sup-convolution} is attained.
The function $U_\lambda$ is continuous, is even in $t$, and is locally Lipschitz
away from $\flatset M$,
with
\[
  0<U_\lambda\le1,\qquad \{U_\lambda=1\}=\flatset M,\qquad U_\lambda(X)\to0\quad \text{as
  } |X|\to\infty .
\]
If $X=(x,t)\notin\flatset M$ then every optimizing pair $(Y,Z)$ has
$U_0(Y)=U_1(Z)=U_\lambda(X)<1$. If in addition $t \neq 0$ then
$\nabla U_0(Y)\ne0$ and $\nabla U_1(Z)\ne0$.
Lastly, given a compact set $E\subset\R^{n+1}\setminus\flatset M$, there is a
constant $C_E$ such that
\begin{equation}\label{eq:convolution-vertical-bound}
  |(U_\lambda)_t(x,t)|\le C_E|t|
\end{equation}
for Lebesgue-almost every $(x,t)\in E$.
\end{lemma}

\begin{proof}
The choice $Y=Z=X$ gives $U_\lambda(X)\ge\min\{U_0(X),U_1(X)\}>0$.
By Lemma~\ref{lem:extended-facts}, $U_0,U_1$ are continuous and vanish
at infinity, so their positive superlevel sets are compact.
A maximizing sequence $(Y_j,Z_j)$ eventually lies in such sets and
therefore has a convergent subsequence, with limit $(Y,Z)$.
The affine relation survives to the limit, and continuity gives
\[
  \min\{U_0(Y),U_1(Z)\}
  =\lim_j\min\{U_0(Y_j),U_1(Z_j)\}=U_\lambda(X).
\]
Thus the supremum in \eqref{eq:sup-convolution} is attained.
The same compactness argument, locally uniformly in $X$, proves upper
semicontinuity of $U_\lambda$. Fixing an optimizing $Y$ and varying $Z$
with $X$ proves lower semicontinuity.
Reflection gives evenness with respect to $t$, and Lemma~\ref{lem:extended-facts}
identifies
$\{U_\lambda=1\}$ with $\flatset M$. Since $\max\{|Y|,|Z|\}\ge|X|$, the decay of each
potential implies
$U_\lambda(X)\to0$.

Suppose an optimizing pair at $X\notin \flatset M$ has
$U_0(Y)<U_1(Z)$.
Then $U_0(Y)<1$ and $Y\notin\flatset K$.
By Lemma~\ref{lem:extended-facts}, a small perturbation of $Y$ raises
$U_0(Y)$. The corresponding admissible perturbation of $Z$ preserves the
strict inequality, contradicting optimality. Interchanging the summands
excludes $U_0(Y)>U_1(Z)$.
Hence $U_0(Y)=U_1(Z)=U_\lambda(X)$. This common value is less than $1$ since $X \notin
\flatset M$, and so $Y\notin \flatset K$ and $Z\notin \flatset L$. Thus both potentials
are smooth near their respective
optimizing inputs.

Suppose in addition that $X$ has nonzero height, $t \neq 0$. If $\grad U_0(Y)=0$ then
formula~\eqref{eq:vertical-derivative} forces $Y$ to have height zero,
and hence $Z$ has nonzero height and $\grad U_1(Z)\ne0$.
By \eqref{eq:horizontal-laplacian}, the Hessian $D^2U_0(Y)$ has a
positive direction. The set of vectors $W\in\R^{n+1}$ with
$D^2U_0(Y)[W,W]>0$ is nonempty, open and symmetric. We may therefore
choose $W$ with $\grad U_1(Z)\cdot W<0$.
For small $h>0$, the pair
\[
  Y+hW,\qquad Z-\frac{1-\lambda}{\lambda}hW,
\]
increases $U_0$ to second order and $U_1$ to first
order while still yielding $X$ as their convex combination, contrary to optimality.
Hence $\grad U_0(Y) \neq 0$, and similarly $\grad U_1(Z) \neq 0$.

On a compact set $E\subset\R^{n+1}\setminus\flatset M$, we have
$U_\lambda\le1-\eta$ for some $\eta>0$. Optimizing inputs $Y,Z$ must stay in some
compact sets
and obey the same bound on their potentials, and hence stay a positive
distance from $\flatset K,\flatset L$. Their potential gradients are
uniformly bounded on the resulting compact sets. Fixing an optimizing $Y$ and replacing
$Z$ by $Z+(X'-X)/\lambda$ gives
\[
  U_\lambda(X')\ge U_\lambda(X)-C|X'-X|
\]
for nearby $X,X'\in E$. The same holds with $X,X'$ interchanged, and so $U_\lambda$ is
Lipschitz
continuous on $E$.

To prove \eqref{eq:convolution-vertical-bound}, write $X=(x,t)\in E$
and take an optimizing pair $Y=(y,t_0)$, $Z=(z,t_1)$. By evenness we may suppose
$t\ge0$. If $t_0<0$ then replacing the heights $t_0,t_1$ by $0,t/\lambda$ decreases the
absolute heights and thus strictly increases both potentials
by \eqref{eq:vertical-derivative}, contradicting optimality of $Y,Z$. Hence $t_0 \geq
0$, and similarly $t_1 \geq 0$.
Hence
\[
  0\le t_0\le\frac{t}{1-\lambda},\qquad
  0\le t_1\le\frac{t}{\lambda}.
\]
The kernel
integrals in \eqref{eq:vertical-derivative} are uniformly bounded along
small vertical shifts of the optimizing inputs. Together with the preceding
estimates on $t_0,t_1$, formula \eqref{eq:vertical-derivative} bounds the vertical
derivative along a shift of magnitude $|h|$ by $C_E(t+|h|)$.
Hence the
change in $U_0,U_1$ is bounded by $C_E(t|h|+h^2)$.
Shifting the heights of $Y,Z$ by $h$ gives an admissible
pair for $(x,t+h)$, and so
\[
  U_\lambda(x,t+h)\ge U_\lambda(x,t)-C_E(t|h|+h^2).
\]
Taking difference quotients with both signs of $h$ implies
\eqref{eq:convolution-vertical-bound} at each point of
differentiability. Those points have full Lebesgue measure by the local Lipschitz
regularity.
\end{proof}

\subsection*{Step 2 --- the local convolution calculation}

The underlying maximum-principle method goes back to
Gabriel \cite[Theorem 1]{G55}.
The gradient and Hessian calculation underlying the following lemma
appears in Borell \cite[proof of Theorem 2.1, pp.~180--181]{B83}.
We use Akman's upper test formulation \cite[Lemma 4.2]{A26}
and reproduce the proof for completeness.
This formulation requires no convexity of level sets.

\begin{lemma}[Aligned gradients and Hessian comparison]\label{lem:local-calculus}
Let $X,Y,Z\in\R^{n+1}$ satisfy $X=(1-\lambda) Y+\lambda Z$,
where $0<\lambda<1$. Suppose
$u,u_0,u_1$ are $C^2$ near $X,Y,Z$ respectively, with $u(X)=u_0(Y)=u_1(Z)$, and $\grad
u_0(Y)$ and $\grad u_1(Z)$ are nonzero, and
\begin{equation}\label{eq:minimality}
  u\bigl((1-\lambda) Y'+\lambda Z'\bigr)
  \ge\min\{u_0(Y'),u_1(Z')\}
\end{equation}
for all nearby $Y',Z'$. Then there is a unit vector $\Xi\in\R^{n+1}$ and
numbers $\ell,\ell_0,\ell_1>0$ with
\begin{equation}\label{eq:gradient-mean}
  \grad u=-\ell\, \Xi, \quad  \grad u_0=-\ell_0\, \Xi, \quad  \grad u_1=-\ell_1\, \Xi,
  \quad \frac1\ell=\frac{1-\lambda}{\ell_0}+\frac{\lambda}{\ell_1}.
\end{equation}
Moreover, as an inequality of symmetric matrices,
\begin{equation}\label{eq:hessian-cubic}
  D^2u\ge
  (1-\lambda)\left(\frac \ell{\ell_0}\right)^{\!\!3}D^2u_0
  +\lambda\left(\frac \ell{\ell_1}\right)^{\!\!3}D^2u_1,
\end{equation}
where the derivatives of $u,u_0,u_1$ are evaluated at $X,Y,Z$, respectively.
\end{lemma}

\begin{proof}
Write $G=\grad u$ and $G_i=\grad u_i \neq 0$, evaluated at $X$ and $Y,Z$ respectively.
First order expansion of \eqref{eq:minimality} in
independent directions $W_0,W_1\in\R^{n+1}$ gives
\[
  G\cdot((1-\lambda)W_0+\lambda W_1)
  \ge\min\{G_0\cdot W_0,G_1\cdot W_1\}.
\]
Set $W_1=0$ and let $W_0$ range over vectors perpendicular to $G_0$.
The inequalities for $W_0$ and $-W_0$ give $G\cdot W_0=0$,
and so $G$ is parallel to $G_0$. Similarly $G$ is parallel to $G_1$.
Write $G=\kappa_i G_i$. Taking \((W_0,W_1)=(G_0,0)\) and \((W_0,W_1)=(0,G_1)\) gives
\(\kappa_0,\kappa_1\ge0\).
Now put $W_i^*=G_i/|G_i|^2$. The inequalities for $\pm(W_0^*,W_1^*)$ give
\[
  (1-\lambda) \kappa_0+\lambda \kappa_1=1 .
\]
Neither $\kappa_0$ nor $\kappa_1$ can vanish, for then $G=0$ would force
both to vanish. Thus $\kappa_i>0$. Taking $\ell=|G|$ and $\ell_i=|G_i|$
gives \eqref{eq:gradient-mean}, with $\Xi=-G/\ell$. Note $\ell=\kappa_i \ell_i$.

It remains to compare the Hessian matrices. Given a vector $W\in\R^{n+1}$, put
\[
  \beta=\frac{\kappa_0^2D^2u_0[W,W]}{2\ell_0},\qquad
  \gamma=\frac{\kappa_1^2D^2u_1[W,W]}{2\ell_1}.
\]
For $h\in\R$, take
\[
  Y_h=Y+h \kappa_0 W+h^2\beta\,\Xi,\qquad
  Z_h=Z+h \kappa_1 W+h^2\gamma\,\Xi.
\]
Their convex combination is
\[
  \begin{split}
    X_h&=(1-\lambda)Y_h+\lambda Z_h\\
    &=X+hW+h^2\bigl((1-\lambda)\beta+\lambda\gamma\bigr) \, \Xi.
  \end{split}
\]
Second order Taylor expansions, using that $\kappa_0G_0=\kappa_1G_1=G$ and
$G_i=-\ell_i\Xi$,
give
\[
  u_0(Y_h)=u(X)+hG\cdot W+o(h^2),\qquad
  u_1(Z_h)=u(X)+hG\cdot W+o(h^2).
\]
Therefore assumption
\eqref{eq:minimality} implies that
\begin{equation}\label{eq:secondorder}
  u(X_h)\ge u(X)+hG\cdot W+o(h^2).
\end{equation}
Expanding the left side of \eqref{eq:secondorder}, canceling the common
zeroth and first order terms and using $G=-\ell \, \Xi$ and
$\kappa_i=\ell/\ell_i$ yields \eqref{eq:hessian-cubic}.
\end{proof}

\subsection*{Step 3 --- the subsolution}

The extended equilibrium
potential $U_K$ satisfies $ \La U_K=0$ by \eqref{eq:weighted-equation}, when $t \neq 0$.
For the supremum--minimum convolution $U_\lambda$ we now prove the subsolution
inequality
\[
  \La U_\lambda=\Delta_{x,t}U_\lambda+\frac a t(U_\lambda)_t\ge0,
  \qquad t\ne0,
\]
first in the viscosity sense and then weakly. After continuation across $t=0$ in
\hyperref[step:crossing-plane]{Step~4} below,
this subsolution inequality will allow the comparison $U_\lambda\le U_M$
in \hyperref[step:comparison-infinity]{Step~5}.

To prove this subsolution inequality, we apply Lemma~\ref{lem:local-calculus} to smooth
upper test functions. Suppose that a $C^2$-smooth function $u$ touches $U_\lambda$ from
above at $X=(x,t)$
with $t>0$. We must show that
\[
  \La u(X)=\Delta_{x,t}u(X)+\frac a t u_t(X)\ge0.
\]
Choose an optimizing pair $Y=(y,t_0)$, $Z=(z,t_1)$, with
\[
  (1-\lambda)t_0+\lambda t_1 = t.
\]
Lemma~\ref{lem:optimizers} supplies
equal values less than $1$, and also nonzero gradients, as needed for applying
Lemma~\ref{lem:local-calculus}. The signs in \eqref{eq:vertical-derivative} and the
alignment
in \eqref{eq:gradient-mean} force $t_0,t_1$ to have the same sign.
They cannot vanish, since $t > 0$. Thus $t_0,t_1>0$ and
$\Xi_{n+1}>0$.

After taking the trace in \eqref{eq:hessian-cubic} and using $\La U_K=\La U_L=0$ to
replace the Laplacians on the right side, one arrives at
\begin{equation}\label{eq:drift-estimate}
  \La u(X)\ge
  a \ell\, \Xi_{n+1} \! \left[
  \ell^2\left(\frac{1-\lambda}{\ell_0^2t_0}
  +\frac{\lambda}{\ell_1^2t_1}\right)-\frac1t\right]
\end{equation}
since $u_t=-\ell \, \Xi_{n+1}$ and $(U_i)_t=-\ell_i \, \Xi_{n+1}$ by Lemma
\ref{lem:local-calculus}.
The bracket on the right side of \eqref{eq:drift-estimate} is nonnegative by
Cauchy--Schwarz, since
\[
  \bigl((1-\lambda) t_0+\lambda t_1\bigr)
  \left(\frac{1-\lambda}{\ell_0^2t_0}+\frac{\lambda}{\ell_1^2t_1}\right)
  \ge\left(\frac{1-\lambda}{\ell_0}+\frac{\lambda}{\ell_1}\right)^{\!2}
  =\frac1{\ell^2}.
\]
Thus $\La u\ge0$, because $a \geq 0$. (This nonnegativity of $a$ is the only
place where $p\ge n-1$ is used in Sections~\ref{sec:extension} and~\ref{sec:upper-proof}.)

We have proved the subsolution inequality
$\La U_\lambda\ge0$ in the viscosity sense, for $t>0$. Reflection gives the same
conclusion for
$t<0$.

We now pass to the weak formulation. By Lemma~\ref{lem:optimizers},
$U_\lambda$ is locally Lipschitz on $\{t\ne0\}$ and hence belongs to
$H^1_{\mathrm{loc}}(\{t\ne0\})$.
On every open set $D\Subset\{t\ne0\}$, the operator $\La$ has smooth
coefficients, the identity principal matrix and no zeroth order term.
Thus Harvey--Lawson \cite[Corollary~9.4]{HL13} implies that
$\La U_\lambda\ge0$ in the sense of distributions.
For nonnegative $\psi\in C_c^\infty(\{t\ne0\})$, testing this distributional
inequality against $|t|^a\psi$ and integrating by parts results in the weak subsolution
inequality
\begin{equation}\label{eq:weak-Ut}
  \int_{\R^{n+1}}|t|^a\grad U_\lambda\cdot\grad\psi\dd X\le0.
\end{equation}

\subsection*{Step 4 --- crossing the extension plane}\label{step:crossing-plane}

We extend \eqref{eq:weak-Ut} across $t=0$ by adapting the zero-flux
reflection argument of Caffarelli--Silvestre \cite[Lemma 4.1]{CS07}
to our subsolution. Let $\psi\ge0$ be smooth with compact support
in $\R^{n+1}\setminus\flatset M$.
Choose an even smooth cutoff $0 \leq \chi_\varepsilon(t) \leq 1$ that equals zero
for $|t|\le\varepsilon$ and $1$ for $|t|\ge2\varepsilon$, with
$|\chi_\varepsilon'|\le C/\varepsilon$. Use
$\psi\chi_\varepsilon$ in \eqref{eq:weak-Ut}.
By the local Lipschitz property in Lemma~\ref{lem:optimizers}, $\grad U_\lambda$ is
essentially bounded on a compact
neighborhood of the support of $\psi$. Since $a>-1$, equivalently $p>n-2$,
dominated convergence handles the term from \eqref{eq:weak-Ut} involving
$\chi_\varepsilon \grad\psi$.
The vertical estimate \eqref{eq:convolution-vertical-bound} bounds the absolute value of
the
term involving $\chi_\varepsilon^\prime \psi$ by
\[
  \frac C\varepsilon\int_{\varepsilon<|t|<2\varepsilon}|t|^{a+1}\dd t
  \le C'\varepsilon^{a+1}\longrightarrow0.
\]
Thus \eqref{eq:weak-Ut} holds throughout
$\R^{n+1}\setminus\flatset M$.

\subsection*{Step 5 --- comparison and the capacity at
infinity}\label{step:comparison-infinity}

We proceed to show $U_\lambda \leq U_M$ and use their asymptotic behavior at infinity to
derive the Brunn--Minkowski inequality.

Lemma~\ref{lem:extended-facts} gives
$\operatorname{div}(|t|^a\grad U_M)=0$ weakly off $\flatset M$.
By Lemmas~\ref{lem:extended-facts} and \ref{lem:optimizers},
$U_M$ and $U_\lambda$ are continuous, locally Lipschitz off $\flatset M$,
equal to $1$ on $\flatset M$ and vanishing at infinity.
Hence, for $\varepsilon>0$,
\[
  w=(U_\lambda-U_M-\varepsilon)_+
\]
is Lipschitz with compact support a positive distance from $\flatset M$.
For sufficiently small $\delta>0$, let $\psi_\delta\ge0$ be a smooth
mollification of $w$, supported in a fixed compact subset of
$\R^{n+1}\setminus\flatset M$.
Testing \eqref{eq:weak-Ut} and the weak equation for $U_M$ with
$\psi_\delta$ and subtracting gives
\[
  \int_{\R^{n+1}}|t|^a\grad(U_\lambda-U_M-\varepsilon)\cdot\grad\psi_\delta\dd X\le0.
\]
The gradients $\grad\psi_\delta$ are uniformly bounded and converge
to $\grad w$ almost everywhere as $\delta\downarrow0$.
Since $|t|^a$ is locally integrable for $a>-1$, that is, for $p>n-2$,
dominated convergence and the chain rule for $w$ give
\[
  \int_{\R^{n+1}}|t|^a|\grad w|^2\dd X\le0.
\]
Hence $\grad w=0$
almost everywhere, and so the compact support forces $w=0$. Letting
$\varepsilon\downarrow0$ gives $U_\lambda \leq U_M$ on $\R^{n+1}$.

Hence by the definition \eqref{eq:sup-convolution} of $U_\lambda$, for all
$Y,Z\in\R^{n+1}$ we have
\begin{equation}\label{eq:compact-potential-comparison}
  U_M((1-\lambda)Y+\lambda Z)
  \ge\min\{U_0(Y),U_1(Z)\}.
\end{equation}
Taking $Y=Z=X$, multiplying by $|X|^p$ and using
\eqref{eq:potential-infinity} yields
\begin{equation}\label{eq:minimumbound}
  \min\{\RCap_p(K),\RCap_p(L)\}
  \leq  \RCap_p(M) .
\end{equation}

For arbitrary nonempty compact $K,L$, put $M=(1-\lambda)K+\lambda L$
and use the parallel neighborhoods $K^r=K+rB_n$, $L^r=L+rB_n$ and
$M^r=M+rB_n$. Applying \eqref{eq:minimumbound} to these neighborhoods
and letting $r\downarrow0$ using Lemma~\ref{lem:outer-continuity}
proves the minimum bound \eqref{eq:minimumbound} for $K,L,M$.

If the capacities of $K$ and $L$ are both positive, Lemma~\ref{lem:outer-continuity} also
allows $r\downarrow0$ in inequality \eqref{eq:compact-potential-comparison}
when $Y,Z$ have positive height. Monotone convergence in the kernel
extends the inequality to nonnegative heights. Replacing the two heights
by their absolute values leaves the right side unchanged and can only
decrease the left side, since the kernel is decreasing with respect to $|t|$. Thus
\eqref{eq:compact-potential-comparison}
holds for all $Y,Z\in\R^{n+1}$.
Taking the supremum over admissible pairs in \eqref{eq:sup-convolution}
gives $U_M\ge U_\lambda$ on $\R^{n+1}$.

The minimum bound \eqref{eq:minimumbound} expresses quasiconcavity under Minkowski
interpolation.
For positive functionals that are homogeneous of degree $1$, quasiconcavity
is equivalent to concavity, as the following standard normalization shows.
Writing $c_0=\RCap_p(K)$ and $c_1=\RCap_p(L)$, if both $c_0$ and $c_1$ are positive then
both
$K/c_0,L/c_1$ have capacity $1$. Applying the minimum bound to these sets
with $\lambda=c_1/(c_0+c_1)$ gives
$\RCap_p((K+L)/(c_0+c_1))\ge1$, and homogeneity yields the Brunn--Minkowski inequality
\eqref{eq:sumBM}.
If either $c_0$ or $c_1$ equals zero, inclusion monotonicity and translation
invariance give \eqref{eq:sumBM}, since $K+L$ contains a translate of
each summand. This proves Theorem~\ref{thm:upper-band} for $p>0$.

\subsection*{Where was \texorpdfstring{$p\ge n-1$}{p >= n-1} used in the proof?}

The hypothesis $p\ge n-1$ is used only to conclude $\La u\ge0$
from \eqref{eq:drift-estimate}. If instead $n-2<p<n-1$ with $n\ge2$,
then $a=p-n+1<0$, and estimate \eqref{eq:drift-estimate} no longer gives the required
sign for a subsolution. The other ingredients in the proof remain valid.
In particular, the cutoff error in
\hyperref[step:crossing-plane]{Step~4} still tends to zero
because $a>-1$.

\section{\bf Proof of Theorem~\ref{thm:upper-equality} ---
equality}\label{sec:upper-equality}

The proof adapts the arguments near infinity of
Caffarelli--Jerison--Lieb \cite[Section 4]{CJL96} and
Akman \cite[Lemma 6.6]{A26} to the weighted extension equation.
The additional work concerns the weight $|t|^a$ and the drift term
$(a/t)U_t$ in the Riesz extension equation.
In \hyperref[step:equality-geometry]{Step~1}, we describe the low level
sets and establish smoothness of the sup-min convolution near infinity.
In \hyperref[step:equality-potentials]{Step~2}, we continue the weak
inequality across $t=0$ and use weighted spherical averages to obtain
equality of potentials near infinity.
In \hyperref[step:equality-rigidity]{Step~3}, we account for the drift term
in the Hessian equality calculation and use $a\ge0$ to deduce homothety.
Analytic continuation and a finite energy measure argument recover the original
sets in \hyperref[step:equality-recovery]{Step~4}, including points outside the supports
of their equilibrium measures.

Suppose $K$ and $L$ are positively homothetic, $L=cK+b$ with $c>0$, and $K$ is convex.
Then
$K+L=(1+c)K+b$, and translation invariance and homogeneity give
equality in \eqref{eq:sumBM}.

For the other direction, first suppose $n=1$ and $p=0$. If equality holds in
\eqref{eq:sumBM} then $K$ and $L$ are positively homothetic and convex (which in one
dimension means they are intervals) by \cite[Theorem 2]{CLline}, or by a special case of
Akman's planar logarithmic equality characterization
\cite[Theorems A(d) and B(e)]{A26}.

\smallskip
From now on, assume $n\ge1$, $p>0$ and $n-1\le p<n$. Suppose equality holds in
\eqref{eq:sumBM}. Normalize as in Section~\ref{sec:upper-proof}: write
$c_0=\RCap_p(K)$ and $c_1=\RCap_p(L)$, replace $K,L$ by
$K/c_0,L/c_1$ and put $\lambda=c_1/(c_0+c_1) \in (0,1)$, retaining the names $K,L$
for the normalized sets. The equality assumption becomes
\[
  \RCap_p(K)=\RCap_p(L)=\RCap_p(M)=1 \quad \text{where} \quad  M=(1-\lambda)K+\lambda L .
\]
Use the extended potentials \eqref{eq:extended-potential}, with
$U_0=U_K$, $U_1=U_L$ and $a=p-n+1\ge0$.
For $E=K,L,M$, the Frostman inequalities and maximum principle used in
Lemma~\ref{lem:extended-facts} give $0<U_E\le1$ everywhere in $\R^{n+1}$ and
$U_E(x,0)=1$ for $\mu_E$-almost every $x\in E$.

\subsection*{Step 1 --- low level sets and smoothness near
infinity}\label{step:equality-geometry}

We describe the low level sets of $U_0,U_1$ and show that $U_\lambda$
is smooth near infinity, preparing the differential comparison in
\hyperref[step:equality-potentials]{Step~2}.

Since the capacities equal $1$, the kernel formula gives
\[
  U_i(X)=|X|^{-p}+O(|X|^{-p-1}),\qquad i=0,1,
\]
as $|X| \to \infty$ with the corresponding differentiated estimates of every order.
Each $U_i$ has a positive lower bound on a fixed large closed ball in $\R^{n+1}$ and is
continuous outside the ball, while
$\partial_rU_i=-pr^{-p-1}+O(r^{-p-2})<0$ outside it, where $r=|X|$.
Hence, for small $s>0$, the level $\{ U_i=s \}$ is a single radial graph.
After rescaling by $s^{1/p}$, this graph converges in $C^2$ to the unit
sphere as $s \to 0$.
Thus $C_i(s)=\{U_i\ge s\}$ is a smooth convex body with positive
principal curvatures, symmetric in $t=0$, and $\grad U_i$ is nonzero
and points inward on its boundary.

We must justify attainment in \eqref{eq:sup-convolution} near infinity. We cannot use
Lemma~\ref{lem:optimizers} because that result requires horizontally regularized sets.
Choose $s_0>0$ so that the preceding level set conclusions hold for
$0<s<s_0$ and $U_0,U_1>s_0$ on the large closed ball.
Taking the supremum in \eqref{eq:compact-potential-comparison}
gives $0<U_\lambda\le U_M$, and so $U_\lambda(X)\to0$
as $|X|\to\infty$. Suppose $0<U_\lambda(X)<s_0$.
For this fixed $X$, any supremizing sequence $(Y_j,Z_j)$ for $U_\lambda(X)$
must be bounded, due to decay of the potentials. Pass to a convergent
subsequence, say $Y_j\to Y$ and $Z_j\to Z$. If $Y$ belongs to the large closed ball then
$U_0(Y) \geq s_0>U_\lambda(X)$, while if $Y$ lies outside the ball then continuity gives
\[
  U_0(Y) = \lim_j U_0(Y_j) \geq
  \lim_j\min\{U_0(Y_j),U_1(Z_j)\}
  =U_\lambda(X).
\]
Thus in both cases $U_0(Y) \geq U_\lambda(X)$. Similarly $U_1(Z) \geq U_\lambda(X)$.
Thus the limiting pair $(Y,Z)$ attains the supremum for $U_\lambda(X)$.

\begin{figure}[!htb]
\centering
\BeginAccSupp{method=pdfstringdef,unicode,Alt={Level set convolution diagram. Three panels show smooth convex sets C0(s), C1(s) and their Minkowski average, the superlevel set of U lambda at level s, with lambda equal to one half. Each set is symmetric about the horizontal extension plane t equals zero, with unequal left and right profiles and no central symmetry. Boundary points Y and Z have the same outward unit normal, denoted by capital Xi. Their midpoint X lies on the boundary of the averaged set and has the same outward normal. The input boundaries are dashed and the averaged boundary is solid. The curves are schematic.},ActualText={Level set convolution diagram. Three panels show smooth convex sets C0(s), C1(s) and their Minkowski average, the superlevel set of U lambda at level s, with lambda equal to one half. Each set is symmetric about the horizontal extension plane t equals zero, with unequal left and right profiles and no central symmetry. Boundary points Y and Z have the same outward unit normal, denoted by capital Xi. Their midpoint X lies on the boundary of the averaged set and has the same outward normal. The input boundaries are dashed and the averaged boundary is solid. The curves are schematic.}}
% Editable figure source: Figures/higher-dim-level-sets.tex
\includegraphics[trim=.5pt .5pt .5pt .5pt]{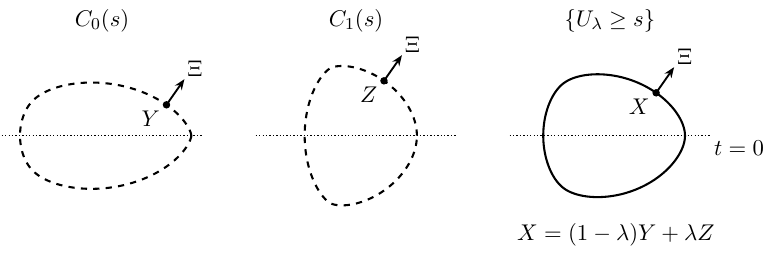}
\EndAccSupp{}
\caption{Schematic of the small-$s$ level sets in $\R^{n+1}$, shown with
$\lambda=1/2$. Boundary points $Y\in\partial C_0(s)$ and
$Z\in\partial C_1(s)$ with common outward normal $\Xi$ combine to
$X=(1-\lambda)Y+\lambda Z \in\partial\{U_\lambda\ge s\}$, with the same normal.
Dotted horizontal lines mark $t=0$, across which the level sets are symmetric.}
\label{fig:level-set-convolution}
\end{figure}

For $0<s<s_0$, attainment when $U_\lambda(X)=s$ and the definition
when $U_\lambda(X)>s$ give
\begin{equation}\label{eq:convolution-level-sets}
  \{U_\lambda\ge s\}=(1-\lambda)C_0(s)+\lambda C_1(s).
\end{equation}
To justify the common normal direction shown in Figure~\ref{fig:level-set-convolution},
consider a unit vector $\Xi$ and let $Y,Z$ be the supporting points of
$C_0(s),C_1(s)$ in direction $\Xi$. Put
\[
  X=(1-\lambda)Y+\lambda Z .
\]
For any $Y'\in C_0(s)$ and $Z'\in C_1(s)$,
\[
  \bigl((1-\lambda)Y'+\lambda Z'\bigr) \cdot \Xi
  \le(1-\lambda) Y \cdot \Xi+\lambda Z \cdot \Xi = X \cdot \Xi.
\]
Thus the hyperplane through $X$ with normal vector $\Xi$
supports the set $\{U_\lambda\ge s\}$.

To show that $U_\lambda$ is smooth near infinity, we use the input
level sets just described. Their support functions $h_i(s,\Xi)$ are
smooth with negative $s$-derivative, as justified in
Appendix~\ref{app:level-set-foliation}. By
\eqref{eq:convolution-level-sets}, the support function of
$\{U_\lambda\ge s\}$ is
\[
  h_\lambda(s,\Xi)=(1-\lambda)h_0(s,\Xi)+\lambda h_1(s,\Xi).
\]

The supporting point of this Minkowski combination is the convex combination of the
input supporting points:
\[
  X(s,\Xi)=(1-\lambda)Y(s,\Xi)+\lambda Z(s,\Xi)
\]
where \(Y,Z\) are the supporting points of \(C_0(s),C_1(s)\) with outward unit normal
\(\Xi\). The derivative of $X(s,\Xi)$ with respect to $\Xi$ on $\Xi^\perp$ is positive
definite, being the weighted
sum of the corresponding positive definite operators for $Y(s,\Xi)$ and $Z(s,\Xi)$,
while its
$s$-derivative has normal component $\partial_s h_\lambda<0$. (For a more detailed
justification of these and the next claims, see Appendix~\ref{app:level-set-foliation}.)
The inverse function theorem therefore shows that these boundaries
form a smooth foliation of a neighborhood of infinity. The level
parameter is $U_\lambda$, proving its smoothness there. Reflection
symmetry gives evenness in $t$.

\subsection*{Step 2 --- equality of potentials near
infinity}\label{step:equality-potentials}

We now prove $U_\lambda=U_M$ near infinity by applying weighted
spherical averages to their nonnegative difference. This equality
will let us deduce homothety of the low level sets in
\hyperref[step:equality-rigidity]{Step~3}.

Put $v=U_M-U_\lambda$. The comparison $U_M\ge U_\lambda$ proved
for arbitrary compact sets in \hyperref[step:comparison-infinity]{Step~5}
of Section~\ref{sec:upper-proof} gives $v\ge0$.
The equality assumption $\RCap_p(K)=\RCap_p(L)=\RCap_p(M)=1$ gives
the same leading term $|X|^{-p}$ for all three potentials.
Choosing $Y=Z=X$ gives $U_\lambda\ge\min\{U_0,U_1\}$, and hence
\[
  0\le v\le U_M-\min\{U_0,U_1\}
  =O(|X|^{-p-1})
\]
uniformly at infinity.

Now we prove $\operatorname{div}(|t|^a\grad v)\le0$ weakly near infinity.
Choose $R_*$ large enough that $U_M$ and $U_\lambda$ are smooth on
$\{|X|>R_*\}$ and the preceding description of the low level sets
applies there.
At a boundary point $X=(x,t)$ of positive height, the corresponding
supporting points $Y\in\partial C_0(s)$ and $Z\in\partial C_1(s)$
have positive heights, by strict convexity and reflection symmetry.
We apply Lemma~\ref{lem:local-calculus} directly with
$u=U_\lambda$, $u_0=U_0$ and $u_1=U_1$ at $X,Y,Z$, respectively.
All three functions are smooth there,
$U_\lambda(X)=U_0(Y)=U_1(Z)=s$, and both $\grad U_0(Y)$ and
$\grad U_1(Z)$ are nonzero.
The local inequality \eqref{eq:minimality} holds by the definition
\eqref{eq:sup-convolution} of $U_\lambda$.
The lemma therefore gives the gradient identities
\eqref{eq:gradient-mean} and Hessian comparison \eqref{eq:hessian-cubic}
pointwise. Taking the trace of the Hessian, then using the partial differential equation
\eqref{eq:weighted-equation} for
$U_0,U_1$ and adding $(a/t)(U_\lambda)_t$ gives
\eqref{eq:drift-estimate} with $u=U_\lambda$.
Its right side is nonnegative, by the same Cauchy--Schwarz calculation
using $a\ge0$, and so $\La U_\lambda\ge0$ for $|X|>R_*$ with $t>0$. Evenness gives the
same inequality for $t<0$.
Since $\La U_M = 0$ by \eqref{eq:weighted-equation}, we have
$\La v\le0$ for $t\ne0$.
Both potentials are smooth and even in $t$, and so $v_t=O(|t|)$ on
compact subsets of $\{|X|>R_*\}$. The cutoff argument in
\hyperref[step:crossing-plane]{Step~4} of Section~\ref{sec:upper-proof}
therefore applies, since $a>-1$, and yields
$\operatorname{div}(|t|^a\grad v)\le0$ weakly throughout $\{|X|>R_*\}$.

Finally, we show that $U_\lambda=U_M$ for $|X|>R_*$. Define a weighted
spherical mean
\[
  m(r)=\int_{\mathbb S^n}|\Theta_{n+1}|^a v(r\Theta)\dd S_\Theta,
  \qquad r>R_*,
\]
where $\Theta$ ranges over the unit sphere $\mathbb S^n\subset\R^{n+1}$.
Clearly $m\ge0$ and $m(r)=O(r^{-p-1})$.
Differentiating under the integral gives
\[
  m'(r)=\int_{\mathbb S^n}|\Theta_{n+1}|^a
  \grad v(r\Theta)\cdot\Theta\dd S_\Theta.
\]
Differentiation is justified here by smoothness of
$v$ and integrability of the angular weight $|\Theta_{n+1}|^a$, since $a>-1$.
For nonnegative $\psi\in C_c^\infty((R_*,\infty))$, use
$\psi(|X|)$ as a test function in the weak inequality for $v$.
Its gradient is $\psi'(r)X/r$, and so polar coordinates give
\[
  0\le\int_{|X|>R_*}|t|^a\grad v\cdot\grad\bigl(\psi(|X|)\bigr)\dd X
  =\int_{R_*}^\infty r^{n+a}m'(r)\psi'(r)\dd r.
\]
Equivalently, $\bigl(r^{n+a}m'(r)\bigr)'\le0$ weakly and hence classically since $m$ is
smooth.
The exponent is $n+a=p+1$, and so the inequality says $r^{p+1}m'(r)$ is decreasing.
Using that $m(s)\to0$ as $s\to\infty$, we find
\[
  -m(r)=\int_r^\infty m'(s)\dd s
  \le r^{p+1}m'(r)\int_r^\infty s^{-p-1}\dd s
  =\frac rp\,m'(r).
\]
Hence $\bigl(r^pm(r)\bigr)'\ge0$ and so $r^pm(r)$ is nonnegative and increasing. But
$r^pm(r)=O(r^{-1})\to0$ and so $m(r)=0$ for all $r>R_*$. Therefore $v=0$ in
$\{|X|>R_*\}$, which means
\begin{equation}\label{eq:equality-convolution}
  U_\lambda=U_M\qquad\text{near infinity in $\R^{n+1}$}.
\end{equation}

\subsection*{Step 3 --- homothety of the low level sets}\label{step:equality-rigidity}

We use equality in the Hessian comparison to show that the low level
sets of $U_0,U_1$ differ by a fixed horizontal translation.
We adapt Colesanti--Salani's equality calculation for radius of curvature
operators \cite[pp.~477--478, (42)--(43)]{CS03} to obtain
\eqref{eq:equality-radius-ratio}, as proved below.

We first prove equality in the Hessian comparison at corresponding
supporting points. Fix $s>0$ small enough that
$\partial\{U_\lambda\ge s\}\subset\{|X|>R_*\}$, where
\eqref{eq:equality-convolution} holds, and choose a unit normal
$\Xi\in\mathbb S^n$ with $\Xi_{n+1}>0$.
Retain the supporting points
$Y=Y(s,\,\Xi)$ of $C_0(s)$ and $Z=Z(s,\,\Xi)$ of $C_1(s)$, and write
\[
  \begin{gathered}
    t_0=Y_{n+1}>0,\qquad t_1=Z_{n+1}>0,\\
    X=(1-\lambda)Y+\lambda Z,\qquad t=(1-\lambda)t_0+\lambda t_1.
  \end{gathered}
\]
Apply Lemma~\ref{lem:local-calculus} to $U_\lambda,U_0,U_1$,
taking all gradients and Hessians at $X,Y,Z$, respectively.
In the lemma's notation,
$\grad U_i=-\ell_i\,\Xi$, $\grad U_\lambda=-\ell\,\Xi$ and
$1/\ell=(1-\lambda)/\ell_0+\lambda/\ell_1$, with $\ell_0,\ell_1,\ell>0$.
Equation~\eqref{eq:hessian-cubic} gives the positive semidefinite matrix
\[
  H=D^2U_\lambda-(1-\lambda)(\ell/\ell_0)^3D^2U_0-\lambda(\ell/\ell_1)^3D^2U_1\ge0.
\]
Applying \eqref{eq:drift-estimate} with $u=U_\lambda$ gives equality
in that estimate, since $\La U_\lambda=\La U_M=0$ by
\eqref{eq:equality-convolution} and its right side is nonnegative.
The estimate was obtained by taking the trace in
\eqref{eq:hessian-cubic}, and so $\operatorname{tr}H=0$.
Since $H$ is positive semidefinite, $H=0$, including when $a=0$.

To express this Hessian equality geometrically, fix $s$ and put
\[
  R_0=D_\Xi Y(s,\,\Xi),\qquad R_1=D_\Xi Z(s,\,\Xi).
\]
These are the positive definite radius of curvature operators on
$\Xi^\perp$ described in Appendix~\ref{app:level-set-foliation}.
By \eqref{eq:convolution-level-sets}, the corresponding operator for
$\partial\{U_\lambda\ge s\}$ is $R=(1-\lambda)R_0+\lambda R_1$.
Differentiating $\grad U_i=-\ell_i\,\Xi$ with respect to $\Xi$
at fixed $s$, the derivative of $\ell_i$ contributes only a normal component.
Taking tangent components gives
$D^2U_i|_{\Xi^\perp}=-\ell_iR_i^{-1}$, and similarly
$D^2U_\lambda|_{\Xi^\perp}=-\ell R^{-1}$.
Here each Hessian has its input restricted to the tangent space \(\Xi^\perp\) and its output orthogonally projected onto that space. Thus $H=0$ gives
\begin{equation} \label{eq:hessianrelation}
  R^{-1}=(1-\lambda)(\ell/\ell_0)^2R_0^{-1}+\lambda(\ell/\ell_1)^2R_1^{-1}.
\end{equation}

A minimum formula now converts this identity into proportionality of the
curvature operators. For each $W\in\Xi^\perp$, we claim that 
\[
  W\cdot R^{-1}W
  =\min_{\substack{W_0,W_1\in\Xi^\perp\\(1-\lambda)W_0+\lambda W_1=W}}
  \bigl[(1-\lambda)W_0\cdot R_0^{-1}W_0
  +\lambda W_1\cdot R_1^{-1}W_1\bigr].
\]
Indeed, the pair $W_i=R_iR^{-1}W$ satisfies the constraint and the first order optimality 
conditions $R_0^{-1}W_0=R_1^{-1}W_1=R^{-1}W$, and so strict convexity makes
it the unique minimizer; substitution then gives the displayed minimum value.
The gradient relation \eqref{eq:gradient-mean} makes $W_i=(\ell/\ell_i)W$ admissible, while the
Hessian identity \eqref{eq:hessianrelation} shows that this pair attains the minimum value $W\cdot R^{-1}W$.
Uniqueness of the minimizer therefore gives $R_iR^{-1}W=(\ell/\ell_i)W$ for all $W$, or
$\ell_iR_i=\ell R$, for $i=0,1$.
In particular, $\ell_0R_0=\ell_1R_1$, and hence
\begin{equation}\label{eq:equality-radius-ratio}
  R_1=cR_0,\qquad c=c(s,\,\Xi)=\ell_0/\ell_1.
\end{equation}

The factor
$c$ is independent of the upper normal direction, as we now prove.
Equation~\eqref{eq:equality-radius-ratio} gives
$\dd Z=c\dd Y$, with $s$ fixed and differentials taken with respect to $\Xi$.
When $n\ge2$, commuting mixed derivatives in local coordinates on the
upper normal hemisphere implies that 
\[
  (\partial_jc)\partial_k Y=(\partial_kc)\partial_j Y.
\]
Since $R_0$ is invertible, the coordinate tangent vectors
$\partial_jY,\partial_kY$ are independent for $j\ne k$, and so $\dd c=0$.
When $n=1$, we have $a=p>0$, and equality in
\eqref{eq:drift-estimate} forces equality in Cauchy--Schwarz, giving
$\ell_0t_0=\ell_1t_1$, hence $t_1=ct_0$.
The last component of $\dd Z=c\dd Y$ gives $\dd t_1=c\dd t_0$.
Differentiating $t_1=ct_0$ therefore yields $t_0\dd c=0$.
Since $t_0>0$, again $\dd c=0$.
Thus for all $n$ we have that $c$ is constant with respect to $\Xi$ on the connected upper normal hemisphere.

Since $c=c(s)$ and $\dd Z=c\dd Y$ by
\eqref{eq:equality-radius-ratio}, we have $\dd(Z-c(s)Y)=0$ at fixed $s$.
Connectedness of the upper normal hemisphere therefore gives a vector
$T(s)$, independent of $\Xi$, such that
\[
  Z(s,\,\Xi)=c(s)Y(s,\,\Xi)+T(s).
\]
As the normal tends to the equator, both supporting points have zero
height by reflection symmetry. Thus $T(s)$ is horizontal, and reflection
gives the same homothety of the full bodies.
The smoothness of $c(s)$ and $T(s)$ follows by fixing
$\Xi=(0,\ldots,0,1)$ and writing
$c(s)=t_1/t_0$ and $T(s)=Z(s,\,\Xi)-c(s)Y(s,\,\Xi)$, since $T(s)$ is horizontal.
The supporting points $Y,Z$ and their heights $t_0,t_1$ are smooth in $s$, with $t_0>0$.

To show that this homothety is independent of $s$, differentiate the
support functions $h_0(s,\Xi)=Y(s,\Xi)\cdot\Xi$ and $h_1(s,\Xi)=Z(s,\Xi)\cdot\Xi$ at fixed $\Xi$. The identities
$U_0(Y(s,\,\Xi))=s$ and $U_1(Z(s,\,\Xi))=s$ give
$\partial_s h_i=-1/\ell_i$. Differentiating
$h_1=c h_0+T\cdot\,\Xi$ and using $c=\ell_0/\ell_1$ gives
\[
  c'(s)h_0(s,\,\Xi)+T'(s)\cdot\,\Xi=0.
\]
At the vertical normal $\Xi=(0,\ldots,0,1)$, we have $h_0=t_0>0$
and $T'\cdot\,\Xi=0$, since $T$ is horizontal. Thus $c'=0$.
The remaining identity $T'\cdot\,\Xi=0$ for every upper unit normal
then gives $T'=0$. Hence both parameters $c$ and $T$ are constant in $s$.

Write $T=(b,0)$ with $b\in\R^n$ and replace $L$ by $L-b$ and $M$ by
$M-\lambda b$, retaining the notation for the sets, their potentials and
supporting points. These translations preserve the Minkowski relation and
\eqref{eq:compact-potential-comparison}. We now have
$Z(s,\,\Xi)=cY(s,\,\Xi)$ and $h_1=ch_0$.

Since $\RCap_p(K)=\RCap_p(L)=1$, both normalized potentials
have leading term $|X|^{-p}$ at infinity, by the kernel formula
\eqref{eq:extended-potential}.
Their low level sets therefore have the same leading radius $s^{-1/p}$,
and so $c=1$. Thus the low level sets coincide, meaning that $U_L=U_K$ near infinity.
Their Minkowski averages are the same sets and so 
\eqref{eq:convolution-level-sets} and \eqref{eq:equality-convolution}
yield that $U_M=U_K$ near infinity.

\subsection*{Step 4 --- recovering the original sets}\label{step:equality-recovery}

Next we show that (after the translation in
\hyperref[step:equality-rigidity]{Step~3}) the normalized sets
$K,L$ coincide and are convex.

The three potentials agree near infinity and are real analytic on the connected upper half-space \(\{t>0\}\), by the kernel formula. Their agreement extends throughout that half-space by the uniqueness theorem for real analytic functions: 
\[
  U_K=U_L=U_M=:U\qquad\text{on }\{t>0\}.
\]
Evenness gives the same identity for $t<0$, and monotone convergence as $t\downarrow0$ extends it to $t=0$.
Put
\[
  F=\{x\in\R^n:U(x,0)=1\},\qquad S=\overline F.
\]
First we identify the common equilibrium support. For $E=K,L,M$,
the equilibrium energy is $V_p(E)=1$, and
\cite[Theorem~4.2.15]{BHS19} gives $U(x,0)\le1$ on $\supp\mu_E$.
The maximum principle \cite[Theorem~4.5.11]{BHS19} therefore gives
$U(x,0)\le1$ throughout $\R^n$.
Together with \cite[Theorem~4.2.14]{BHS19}, this bound gives $U(x,0)=1$
on $E$ except possibly on a set of zero $p$-capacity.
Finite energy measures do not charge such sets \cite[Lemma~4.2.13]{BHS19},
and so $\mu_E(F)=1$. Every neighborhood of a point of $\supp\mu_E$
therefore meets $F$, giving $\supp\mu_E\subset S$.
Conversely, since $p>0$ and $p>n-2$, the function $U(\cdot,0)$ is
strictly subharmonic on $\R^n \setminus \supp\mu_E$
\cite[proof of Theorem~4.5.10(ii)]{BHS19}. 
Since $U\le1$, the strong maximum principle implies $U<1$ outside $\supp\mu_E$. Thus $F\subset\supp\mu_E$,
and so $S\subset\supp\mu_E$ by closedness. Hence
\[
  S=\supp\mu_K=\supp\mu_L=\supp\mu_M\subset K\cap L\cap M.
\]
Each measure is an equilibrium measure on $S$ with energy $1$, since
restricting the admissible measures to $S$ preserves its minimality.
Uniqueness of equilibrium measures \cite[Theorem~4.4.5]{BHS19} forces
$\mu_K=\mu_L=\mu_M=: \mu$.

Akman \cite[proof of Proposition 6.17]{A26} also uses dilations toward a point
of one summand in recovering the filled sets. Here a finite energy pushforward
measure yields the required support inclusion directly.
To recover every point of $K$, fix $x\in K$ and let $\nu_x$ be the
pushforward of $\mu$ under $y\mapsto(1-\lambda)x+\lambda y$.
Since $S\subset L$, this probability measure is supported on $M$, and scaling gives
\[
  \iint_{M\times M}|z-z'|^{-p}\dd\nu_x(z)\dd\nu_x(z')
  =\lambda^{-p}<\infty, 
\]
so that $\nu_x$ has
finite energy. The exceptional set for $U_M=1$ has zero $p$-capacity, and thus $\nu_x(F)=1$. Hence
\[
  (1-\lambda)x+\lambda S=\supp\nu_x\subset S.
\]
Choose $z\in S$ and iterate this inclusion to obtain
$(1-\lambda^j)x+\lambda^jz\in S$ for every $j\ge1$.
These points tend to $x$, and closedness gives $x\in S$.
Thus $K=S$. Similarly $L=S$.

The same inclusion now proves convexity. Since $K=S$, the set $S$ is
closed under the fixed average $(1-\lambda)y+\lambda z$ for $y,z\in S$.
The fixed-average argument in Akman \cite[proof of Lemma 6.7]{A26}
shows that the closed set $S$ is convex: a gap in its intersection with
a segment joining two points of $S$ would have endpoints in $S$,
whose fixed average would lie inside the gap. Hence $K=L$ is convex. 
Undoing normalization and translation proves convexity and positive
homothety of the original sets.

\subsection*{Where did the equality proof use \texorpdfstring{$p\ge n-1$}{p >= n-1}?}

The condition $p\ge n-1$ is used in
\hyperref[step:equality-potentials]{Step~2} and
\hyperref[step:equality-rigidity]{Step~3} through the sign $a\ge0$.
The cutoff and spherical integral arguments in
\hyperref[step:equality-potentials]{Step~2} require only $p>n-2$
within the standing range $0<p<n$.
Once the potential identities are available,
the recovery in \hyperref[step:equality-recovery]{Step~4} also requires only
$p>n-2$.
When $n\ge2$, the arguments work for $p \geq n-1$. When $n=1$, the height argument uses the strict assumption that $a=p>0$.

\subsection*{Proof of Corollary~\ref{cor:width-diameter}}

The inequalities for $n=1$ were explained after the corollary.
For $n\ge2$, we follow Borell's rotation averaging argument
\cite[Examples 7.4 and 7.5]{B84}. Put $C=\operatorname{conv}K$.
Hadwiger's rotation means theorem
\cite[Theorem 3.3.5]{S14} supplies rotations $Q_{j,k}$ of $\R^n$ whose
Minkowski averages
\[
  A_j=\frac1{m_j}\bigl(Q_{j,1}C+\cdots+Q_{j,m_j}C\bigr)
\]
converge in Hausdorff distance to a ball as $j \to \infty$.
Additivity and rotational invariance of mean width give $w(A_j)=w(K)$,
and so the limiting ball has radius $r=w(K)/2$.
Our Brunn--Minkowski Theorem \ref{thm:upper-band} and orthogonal invariance of capacity
give
$\RCap_p(A_j)\ge\RCap_p(C)\ge\RCap_p(K)$.
Given $\varepsilon>0$, Hausdorff convergence places
$A_j$ in a ball of radius $r+\varepsilon$ for all large $j$, so that
$\RCap_p(K)\le(r+\varepsilon)\RCap_p(B_n)$.
Letting $\varepsilon\downarrow0$ proves the mean width inequality in the corollary.
The mean width is at most the diameter and they are equal precisely when $C$ has
constant width. Thus the inequalities in the corollary are proved.

Balls attain both bounds. Conversely, in any dimension $n\ge1$,
suppose equality holds in the mean width bound.
Since $w(K)>0$, we have $\RCap_p(K)>0$.
For an orthogonal map $Q$, put $H=(K+QK)/2$.
Brunn--Minkowski and the bound just proved imply
\[
  \RCap_p(K)\le\RCap_p(H)
  \le\frac{w(H)}2\RCap_p(B_n)=\RCap_p(K),
\]
because $w(H)=w(K)$.
Equality throughout and Theorem~\ref{thm:upper-equality} force $K$
to be convex and $QK$ to be a translate of $K$, since their capacities
agree.
Center the smallest enclosing ball of \(K\) at the origin. Writing \(QK=K+b\) for the
translate of $K$, its smallest enclosing ball is centered at the origin by
orthogonality of $Q$ and also at \(b\) by translation. Uniqueness of the smallest ball
therefore requires \(b=0\) and so \(QK=K\), for every orthogonal \(Q\). Convexity then
makes $K$ a ball.

Equality between the first and last quantities in
\eqref{eq:width-diameter} would imply equality in the first, and so again $K$ would be a
ball.

\section{\bf Planar counterexamples for negative
exponents}\label{sec:negative-planar-proof}

We first use the vertices of a square to treat $p\le-1$, and then explain how known
isodiametric counterexamples rule out a further range of negative exponents.

\subsection*{A square counterexample}

\begin{example}[Perpendicular two point summands]\label{ex:square}
If $p\le-1$, then \eqref{eq:sumBM} fails for $K=\{0,e_1\}$ and
$L=\{0,e_2\}$, where $e_1,e_2$ are the coordinate unit vectors in $\R^2$.
\end{example}

\begin{proof}
Put $s=-p\ge1$. Then $\RCap_p(K)=\RCap_p(L)=2^{-1/s}$, by balancing the masses at the
two points. The sum $K+L$ consists of the four vertices of a unit square. Suppose the
total mass on one diagonal is
$\tau\in[0,1]$. The interaction energy between diagonals
is $2\tau(1-\tau)$, since vertices on different diagonals are distance one apart.
The contribution within each diagonal is proportional to the product of
its endpoint masses, and so is maximized by equal division.
The resulting measure $\mu_\tau$ assigns mass $\tau/2$ to each endpoint
of one diagonal and $(1-\tau)/2$ to each endpoint of the other, giving energy
\[
  \iint |x-y|^s\dd\mu_\tau(x)\dd\mu_\tau(y)
  =2\tau(1-\tau)+2^{s/2-1}\tau^2+2^{s/2-1}(1-\tau)^2.
\]
Maximizing over $\tau\in[0,1]$ gives
\[
  \begin{split}
    \RCap_p(K+L)^s
    &=\begin{cases}
      2^{s/2-2}+2^{-1},&1\le s< 2,\\
      2^{s/2-1},&s \geq 2,
    \end{cases}
    \\
    & <2^{s-1}
    =\bigl(\RCap_p(K)+\RCap_p(L)\bigr)^{\!s} ,
  \end{split}
\]
since when $1 \leq s < 2 $ the expression is less than $1$ and hence less than
$2^{s-1}$, and when $s \geq 2$ one uses $s/2 < s$. Thus $ \RCap_p(K+L) <
\RCap_p(K)+\RCap_p(L)$ and so the Brunn--Minkowski inequality fails when $p \leq -1$.
\end{proof}

\subsection*{The isodiametric obstruction}

The next corollary converts known counterexamples to capacitary
isodiametry into convex counterexamples to Brunn--Minkowski.

\begin{corollary}[Failure from the isodiametric
inequality]\label{cor:isodiametric-failure}
Let $q_*=-0.856\ldots$ be the unique root in $(-1,0)$ of
\[
  2\sqrt\pi\,\Gamma(1-q_*/2)
  =3\Gamma((1-q_*)/2).
\]
If $q<q_*$, then \eqref{eq:sumBM} (with $q$ instead of $p$) fails for some nonempty
convex compact sets $K,L\subset\R^2$.
\end{corollary}

\begin{proof}
Let $D=\frac12B_2$ be the closed disk of diameter $1$, and let $T$ be
the vertex set of an equilateral triangle with side length one.
This $T$ has capacity $(2/3)^{-1/q}$; see
\cite[Section 6]{CLratios}. For $-2<q<q_*$, this capacity exceeds $\RCap_q(D)$
by taking $p\to-\infty$ in the capacity-ratio calculations of
\cite[Theorem 7 and its proof]{CLratios} so that the denominator
$\RCap_p$ tends to diameter; see also \cite[Section 5]{F26}.
For $q\le-2$, we have directly that $\RCap_q(D)=2^{1/q}<(2/3)^{-1/q}$
by \cite[Appendix A]{CLproperties}.

Passing to the convex triangle $K=\operatorname{conv}T$ preserves diameter
and can only increase capacity. Thus $\diam K=1$ and
\[
  \RCap_q(K)\ge\RCap_q(T)>\RCap_q(D).
\]
Take $L=-K$. The difference of two points in $K$ has magnitude at most $1$,
and so $K+L\subset B_2$. Hence
\[
  \begin{aligned}
    \RCap_q(K+L)
    & \le\RCap_q(B_2)\\
    &=2\RCap_q(D) \\
    & <2\RCap_q(K) =\RCap_q(K)+\RCap_q(L).
  \end{aligned}
\]
This proves failure of Brunn--Minkowski for the convex compact summands
$K,L$.
\end{proof}

\begin{problem}[Planar negative exponents]\label{prob:planar-negative}
Does there exist an exponent $q<0$ for which the Brunn--Minkowski inequality
\eqref{eq:sumBM} (with $p=q$) holds
for every pair of nonempty compact sets $K,L\subset\R^2$?
Does the situation change if convexity is assumed?
\end{problem}

\section{\bf Supplementary bounds from couplings}\label{sec:transport}

Optimal transport with cost functional $|x-y|^2$ gives an $\ell^2$-type lower
bound for capacity. The bound follows from Brunn--Minkowski when
that inequality holds, but the independent proof below also covers a range of exponents
$p \geq -2$
for which Brunn--Minkowski is either open or false.
In addition, for $p<-2$ the theorem includes a weaker bound that is obtained through a
simple coupling.

\begin{theorem}[Lower bounds from couplings]\label{thm:transport}
Let $n\ge1$ and suppose $K,L\subset\R^n$ are nonempty compact sets.
If $-2\le p<n$, then
\begin{align}
  \RCap_p(K+L)&\ge
  \bigl(\RCap_p(K)^2+\RCap_p(L)^2\bigr)^{1/2}.
  \label{eq:transport-bound}\\
  \intertext{If $p<-2$, then}
  \RCap_p(K+L)&\ge
  \bigl(\RCap_p(K)^{-p}+\RCap_p(L)^{-p}\bigr)^{\!-1/p}.
  \label{eq:product-bound}
\end{align}
\end{theorem}

The investigation of equality cases is left to the interested reader.

For $n=1$, the theorem follows from the
Brunn--Minkowski inequality in Theorem~\ref{thm:line}.

\begin{proof}
If either capacity vanishes, the asserted bound follows because
$K+L$ contains a translate of each summand. Otherwise put
$c_0=\RCap_p(K)>0$ and $c_1=\RCap_p(L)>0$.
Let $\mu,\nu$ be equilibrium probability measures on $K,L$, respectively.
Then
\[
  \begin{aligned}
    \iint_{K\times K}|x-x'|^{-p}\dd\mu(x)\dd\mu(x')&=c_0^{-p},\\
    \iint_{L\times L}|y-y'|^{-p}\dd\nu(y)\dd\nu(y')&=c_1^{-p},
  \end{aligned}
\]
for $p\ne0$. For $p=0$, the corresponding logarithmic energies equal
$\log(1/c_0)$ and $\log(1/c_1)$.

A coupling $\pi$ of $\mu,\nu$ is a probability measure on $K\times L$
with marginals $\mu,\nu$ \cite[Definition 1.1]{V09}. Choose a coupling that minimizes
the cost
$\int_{K\times L}|x-y|^2\,d\pi(x,y)$
\cite[Theorem 4.1]{V09}.
The cost is continuous, nonnegative and bounded on $K\times L$, and so
\cite[Theorem 5.10(ii)]{V09} supplies a closed set
$\Gamma\subset K\times L$ with $\pi(\Gamma)=1$ that is cyclically
monotone for this cost. Since $\Gamma$ is closed and has full
$\pi$-measure, it contains $\supp\pi$. Taking a cycle of length two in
\cite[Definition 5.1]{V09} therefore gives
\[
  |x-y|^2+|x'-y'|^2\le |x-y'|^2+|x'-y|^2
\]
for all $(x,y),(x',y')\in\supp\pi$. Expanding and cancelling yields
\begin{equation}\label{eq:cyclical}
  (x-x')\cdot(y-y')\ge0.
\end{equation}

First suppose $-2\le p<n$ and set
\[
  c=\sqrt{c_0^2+c_1^2},\qquad \lambda=\frac{c_1^2}{c^2}.
\]
On $(\supp\pi)^2$, define the nonnegative functions
$D_0=|x-x'|/c_0$, $D_1=|y-y'|/c_1$ and
$D=|(x+y)-(x'+y')|/c$.
The inequality \eqref{eq:cyclical} implies that
\[
  D^2\ge(1-\lambda)D_0^2+\lambda D_1^2.
\]
For $p>0$, the function $r\mapsto r^{-p/2}$ is decreasing and convex,
while for $-2\le p<0$ it is increasing and concave (linear at $p=-2$);
this concavity is where we use $p\ge-2$. These properties along with the preceding
displayed equation yield that
\[
  D^{-p}
  \begin{cases}
    \le(1-\lambda) D_0^{-p}+\lambda D_1^{-p},&p>0,\\
    \ge(1-\lambda) D_0^{-p}+\lambda D_1^{-p},&-2\le p<0.
  \end{cases}
\]
Integrate these inequalities with respect to $\pi\otimes\pi$ and note that each
$D_j^{-p}$ has integral $1$ for $j=0,1$.

Now let $\eta$ be the pushforward of $\pi$ under $(x,y)\mapsto x+y$,
so that $\eta\in\PP(K+L)$. (It is the midpoint interpolant in McCann's
displacement interpolation \cite[Definition 1.1 and Remark 1.8]{McC97},
dilated by $2$.)
The energy of $\eta$ is
\[
  \int_{K+L}
  \int_{K+L}
  |z-z'|^{-p}\,d\eta(z)d\eta(z')
  =
  c^{-p}
  \int_{K\times L}
  \int_{K\times L}
  D^{-p}\,d(\pi\otimes\pi) ,
\]
by writing $z=x+y,z'=x'+y'$. This last expression is at most $c^{-p}$ for
$p>0$ and at least $c^{-p}$ for $-2\le p<0$.
Raising this energy estimate to the power $-1/p$ yields $\RCap_p(K+L)\ge c$ in both
cases, proving \eqref{eq:transport-bound}.

For $p=0$, concavity of the logarithm applied to $D^2$ gives
\[
  \log D
  \ge(1-\lambda)\log D_0+\lambda\log D_1.
\]
Both terms on the right have integral zero. They are integrable because
the logarithmic equilibrium energies are finite and the sets are compact.
The trial measure $\eta$ consequently has logarithmic energy at most $\log 1/c$, giving
logarithmic capacity at least $c$.

Finally suppose $p<-2$ and write $s=-p>2$.
Use the independent coupling $\pi=\mu\otimes\nu$ on $K \times L$ and write
$a=x-x'$ and $b=y-y'$.
Interchanging $y$ and $y'$ and then applying convexity gives
\begin{align*}
  \int_{(K\times L)^2}|a+b|^s\dd(\pi\otimes\pi)
  &=\frac12\int_{(K\times L)^2}
  \bigl(|a+b|^s+|a-b|^s\bigr)\dd(\pi\otimes\pi)\\
  &\ge\int_{(K\times L)^2}
  \bigl(|a|^2+|b|^2\bigr)^{\!s/2}\dd(\pi\otimes\pi)\\
  &\ge c_0^s+c_1^s,
\end{align*}
where the last inequality uses $s/2>1$ and the equilibrium energies of
$\mu,\nu$. Again let $\eta$ be the pushforward of $\pi$ under addition.
Then $\eta\in\PP(K+L)$, and so its energy is at most $V_p(K+L)$,
proving \eqref{eq:product-bound}.
\end{proof}

\appendix
\section{\bf Supporting points and smooth level sets}\label{app:level-set-foliation}

We expand the supporting point calculation in
\hyperref[step:equality-geometry]{Step~1} of
Section~\ref{sec:upper-equality}, explaining why the low level sets of
$U_\lambda$ give smooth coordinates near infinity.
The construction follows Caffarelli--Jerison--Lieb \cite{CJL96} and
Akman \cite{A26}, with support function coordinates as in
Colesanti--Salani \cite{CS03}. Section~\ref{sec:upper-equality} gives
more detailed credits.
We use only the smoothness, positive curvature and
nonvanishing gradients of the low level sets of $U_0,U_1$ that are established there.

For $0<s<s_0$, retain the notation $C_i(s)=\{U_i\ge s\}$, $i=0,1$,
and fix $0<\lambda<1$. For $\Xi\in\mathbb S^n$, let $Y(s,\Xi)$ and
$Z(s,\Xi)$ be the unique supporting points of $C_0(s),C_1(s)$ with
outward unit normal $\Xi$. Write
\[
  C_\lambda(s)=(1-\lambda)C_0(s)+\lambda C_1(s),\qquad
  X(s,\Xi)=(1-\lambda)Y(s,\Xi)+\lambda Z(s,\Xi).
\]
The supporting point argument in Section~\ref{sec:upper-equality} shows
that $X(s,\Xi)$ is the unique supporting point of $C_\lambda(s)$ in
direction $\Xi$. Its support function is
\[
  h_\lambda(s,\Xi)=X(s,\Xi)\cdot\Xi
  =(1-\lambda)h_0(s,\Xi)+\lambda h_1(s,\Xi).
\]

\subsection*{Varying the normal direction}

Fix $s$. Since $\Xi$ is a unit vector, its infinitesimal variations
belong to $\Xi^\perp$, the tangent space of $\mathbb S^n$ at $\Xi$.
The resulting variations of $Y,Z$ are tangent to their boundaries,
and so belong to the same vector space $\Xi^\perp$.
For fixed $s$, the radius of curvature operators are
\[
  R_0=D_\Xi Y(s,\Xi):\Xi^\perp\to\Xi^\perp,\qquad
  R_1=D_\Xi Z(s,\Xi):\Xi^\perp\to\Xi^\perp
\]
used in \hyperref[step:equality-rigidity]{Step~3}.
We explain their positivity by a local graph calculation; see also
\cite[Section 2.5, pp.~113--117]{S14} for the standard terminology.

Near $Y \in \partial C_0(s)$, express $\partial C_0(s)$ as a graph $W \mapsto
Y+W-f(W)\Xi$ where $W\in\Xi^\perp$ and $f(0)=0, \grad f(0)=0$. The outward unit normal
of this graph is
\[
  N(W)=\frac{\Xi+\grad f(W)}{\sqrt{1+|\grad f(W)|^2}},
\]
where $\grad f(W)\in\Xi^\perp$. Thus $DN(0)=D^2f(0)$.
The eigenvalues of this symmetric matrix are the positive principal
curvatures. The inverse function theorem therefore lets us recover
the boundary point smoothly from its normal. This recovery first maps a
normal direction to its graph coordinate using $N^{-1}$, then applies
the graph parametrization $W\mapsto Y+W-f(W)\Xi$.
At $W=0$, the derivative of this parametrization is the identity on
$\Xi^\perp$, since $\grad f(0)=0$. Thus the chain rule and the derivative
of an inverse map give
\[
  R_0=I\,D(N^{-1})(\Xi)
  =I\,(DN(0))^{-1}=(D^2f(0))^{-1}.
\]
The same calculation applies at $Z$.
Consequently $R_0,R_1$ are symmetric positive definite operators,
with eigenvalues equal to the principal radii of curvature.
Differentiating the formula for $X$ yields
\begin{equation}\label{eq:appendix-normal-derivative}
  D_\Xi X=R:=(1-\lambda)R_0+\lambda R_1.
\end{equation}
As a convex combination of positive definite operators, \(R\) is positive definite and
hence invertible on \(\Xi^\perp\).

\subsection*{Varying the level}

Since $\grad U_i\ne0$, the level graphs depend smoothly on $s$.
Positive curvature makes the normal map's derivative invertible, giving
joint smoothness of $Y,Z$ in $(s,\Xi)$ by the inverse function theorem
with parameter $s$; see also Colesanti--Salani \cite[Theorem 4]{CS03}.

Now hold $\Xi$ fixed and vary $s$. Put
\[
  \ell_0=|\grad U_0(Y)|>0,\qquad
  \ell_1=|\grad U_1(Z)|>0.
\]
The gradients point inward, and so
$\grad U_0(Y)=-\ell_0\Xi$ and $\grad U_1(Z)=-\ell_1\Xi$.
Differentiating $U_0(Y(s,\Xi))=s$ gives
\[
  -\ell_0\Xi\cdot\partial_sY=1.
\]
Since $h_0(s,\Xi)=Y(s,\Xi)\cdot\Xi$, we obtain
$\partial_sh_0=-1/\ell_0$, and similarly $\partial_sh_1=-1/\ell_1$.
Therefore
\begin{equation}\label{eq:appendix-level-velocity}
  \Xi\cdot\partial_sX=\partial_sh_\lambda
  =-\frac{1-\lambda}{\ell_0}-\frac{\lambda}{\ell_1}<0.
\end{equation}
This equality describes the normal component of the velocity
$\partial_sX$: increasing the level moves the supporting point inward
in its normal component. The velocity may also have a tangential
component.

\subsection*{The inverse function theorem and the foliation}

For an infinitesimal variation $(\dot s,W)\in\R\times\Xi^\perp$,
the full derivative of $(s,\Xi)\mapsto X(s,\Xi)$, using
\eqref{eq:appendix-normal-derivative}, is
\[
  (\dot s,W)\mapsto\dot s\,\partial_sX+RW.
\]
If this vector vanishes then taking the scalar product with \(\Xi\) and recalling that
\(RW\in\Xi^\perp\) gives \(\dot s\,\partial_s h_\lambda=0\), and hence $\dot s=0$ by
\eqref{eq:appendix-level-velocity}. Then $RW=0$ forces $W=0$ by positive definiteness.
The derivative is therefore
invertible. Using local coordinates for \(\Xi\) on the unit sphere \(\mathbb S^n\), the
inverse function theorem gives a smooth local inverse whose coordinates are the level
\(s\) and the normal direction \(\Xi\).

We verify that these coordinates cover a neighborhood of infinity.
Integrating the negative derivative in
\eqref{eq:appendix-level-velocity} gives a continuous
gap $h_\lambda(s,\Xi)-h_\lambda(s',\Xi)>0$ when $s<s^\prime$. The gap has a positive
minimum for $\Xi \in {\mathbb S^n}$, and so
\[
  C_\lambda(s')\subset\operatorname{int}C_\lambda(s)
  \qquad\text{if }0<s<s'<s_0.
\]
If $s_j\to s_*\in(0,s_0)$ and $P\in C_\lambda(s_j)$ for every $j$, then
$P\cdot\Xi\le h_\lambda(s_j,\Xi)\to h_\lambda(s_*,\Xi)$ for every unit $\Xi$,
and hence $P\in C_\lambda(s_*)$.
If $P\in\operatorname{int}C_\lambda(s)$, the uniform positive margin in
these inequalities persists for nearby levels by continuity. The bodies contain balls of
radii
tending to infinity as $s\downarrow0$, by the asymptotics of $U_0,U_1$.
Fix $0<s_1<s_0$. For $P$ outside $C_\lambda(s_1)$, take the largest
level $s_*$ for which it belongs to $C_\lambda(s)$.
The limit property gives $P\in C_\lambda(s_*)$ with $s_*<s_1$.
An interior point would remain in the body at slightly larger levels,
contradicting maximality; hence $P\in\partial C_\lambda(s_*)$.
Strict nesting makes this boundary level unique.
Thus the boundaries with $0<s<s_1$ cover
$\R^{n+1}\setminus C_\lambda(s_1)$ without intersecting one another.
In the smooth local coordinates just proved, each boundary is a slice
with $s$ constant. This is the asserted smooth foliation of a
neighborhood of infinity. In particular, the foliation gives the level as a smooth
function \(s=s(X)\) of the point $X$.

Finally, \eqref{eq:convolution-level-sets} gives
$C_\lambda(s)=\{U_\lambda\ge s\}$. At a boundary point $X \in \partial C_\lambda(s)$ we
have $U_\lambda(X)\ge s$. A larger value would place $X$ in the interior of
$C_\lambda(s)$, a contradiction, and so $U_\lambda(X) = s = s(X)$. Hence \(U_\lambda\)
is smooth in a neighborhood of infinity.

Reflection in $t=0$ preserves every $C_\lambda(s)$,
and so this level parameter is even in $t$.

\section*{\bf Acknowledgments}

Richard Laugesen's research was supported by grants from the Simons
Foundation (\#964018) and the National Science Foundation (\#2246537),
and by the University of Illinois Campus Research Board award RB26003.

This work was supported in part by a grant of access to OpenAI models
through the ChatGPT for Academic Researchers program.

We thank Almut Burchard for suggesting in 2024 an optimal transport approach to
capacitary Brunn--Minkowski inequalities, which we have pursued in
Section~\ref{sec:transport}.

\section*{\bf AI usage statement}

ChatGPT 6.1 helped develop the mathematical arguments and carry
out numerical checks and literature searches. It also assisted with
drafting and revising the manuscript.

\raggedbottom

\end{document}